\documentclass[12pt]{amsart}
\usepackage{hyperref}

\usepackage{amssymb}
\usepackage{graphicx}
\usepackage{amscd}
\usepackage{url}
\usepackage[utf8]{inputenc}
\usepackage[T1]{fontenc}

\theoremstyle{plain}
\newtheorem{theorem}{Theorem}
\newtheorem*{acknowledgement}{Acknowledgement}

\newtheorem{conjecture}[theorem]{Conjecture}
\newtheorem{corollary}[theorem]{Corollary}

\newtheorem{lemma}[theorem]{Lemma}
\newtheorem*{notation}{Notation}

\newtheorem{proposition}[theorem]{Proposition}
\theoremstyle{definition}
\newtheorem{remark}[theorem]{Remark}
\newtheorem*{question}{Question}
\newtheorem{definition}[theorem]{Definition}

\newtheorem{example}[theorem]{Example}

\numberwithin{equation}{section}
\numberwithin{theorem}{section}
\newcommand{\codim}{\text{codim}}

\def\C{{\mathbb C}}
\newcommand{\Z}{{\mathbb Z}}
\newcommand{\R}{{\mathbb R}}
\newcommand{\Q}{{\mathbb Q}}
\def\ol#1{\overline{#1}}
\def\tR{\widetilde{L}}
\def\p{\partial}

\newcommand{\PV}{\mathbb{P}V}
\def\PP{\mathbb{P}}

\DeclareMathOperator{\ord}{ord}
\DeclareMathOperator{\MVol}{MVol}

\title[Puiseux coefficients]{Puiseux coefficients and parametric deformations of plane curve singularities}
\author{Maciej Borodzik}
\address{Institute of Mathematics, University of Warsaw, ul. Banacha 2,
02-097 Warsaw, Poland}
\email{mcboro@mimuw.edu.pl}
\date{\today}
\makeatletter
\@namedef{subjclassname@2010}{\textup{2010} Mathematics Subject Classification}
\makeatother
\subjclass[2010]{Primary 14B07; Secondary 14H20, 14H50, 14H15}
\keywords{plane curve singularity, parametric deformation, Puiseux expansion}
\thanks{Supported by Polish  MNiSz Grant N N201 397937. The author is also supported by Foundation for Polish Science (FNP)}
\begin{document}
\begin{abstract}
We study deformations of cuspidal singularities of plane curves by looking at the coefficients of the Puiseux expansions regarded
as functions of the deformation parameter. 
We obtain new concrete results concerning the adjacency problem.
Moreover, we show some rather unexpected properties of Puiseux coefficients treated as functions on a suitably defined parameter space.
The methods used in paper are elementary.
\end{abstract}
\maketitle

\section{Introduction}
\subsection{Preliminaries}\label{sec:prel}
A cuspidal plane curve singularity is a germ of a plane curve in $\C^2$, locally parametrized as $C=\{(x(t),y(t))\in\C^2\colon t\in D\}$, where
$D\subset\C$ is a disk with center at zero and $x(t)$ and $y(t)$ are analytic functions with the following expansion in power series.

\begin{equation}\label{eq:param}
\begin{split}
x(t)&=a_pt^p+a_{p+1}t^{p+1}+\dots\\
y(t)&=b_qt^q+b_{q+1}t^{q+1}+\dots.
\end{split}
\end{equation}

Here $a_p,a_{p+1},\dots$ and $b_q,b_{q+1},\dots$ are complex coefficients.
Throughout the paper we shall assume that the singularity occurs at $(0,0)\in\C^2$, and, unless specified otherwise, that $a_pb_q\neq 0$, 
so that $p$ is the order of $x(t)$ at $0$
and $q$ is the order of $y$. 
To study the singularity at the origin, it is convenient to consider the following Puiseux expansion

\begin{equation}\label{eq:puis}
y=c_qx^{q/p}+c_{q+1}x^{(q+1)/p}+\dots.
\end{equation}

The Puiseux coefficients $c_q,c_{q+1},\dots$ can be expressed as explicit functions of parameters $a_p,a_{p+1},\dots$ and $b_q,b_{q+1},\dots$. More
precisely, for any $i\ge 0$, we can write

\begin{equation}\label{eq:cqi}
c_{q+i}=a_p^{-(q+i)/p-i}\gamma_{q+i}
\end{equation}
and then $\gamma_{q+i}$ is a polynomial in the coefficients $a_p,\dots$ and $b_q,\dots$. For example

\begin{equation}\label{eq:stupidgamma}
\begin{split}
\gamma_{q}&=b_q\\
\gamma_{q+1}&=a_pb_{q+1}-\frac{q}{p}a_{p+1}b_q\\
\gamma_{q+2}&=a_p^2b_{q+2}-\frac{q+1}{p}a_pa_{p+1}b_{q+1}+\frac{q^2+qp+2q}{2p^2}a_{p+1}^2b_q-\frac{q}{p}a_pa_{p+2}b_q.
\end{split}
\end{equation}

In general, $\gamma_{q+i}$ is a sum of monomials of the form $\eta_{k_p\ldots k_{p+i}}\cdot a_p^{k_p}a_{p+1}^{k_{p+1}}\dots a_{p+i}^{k_{p+i}}b_j$, where $k_p+\dots+k_{p+i}=i$,
$j+k_{p+1}+2k_{p+2}+\dots+ik_{p+i}=i$ and $\eta_{k_p\ldots k_{p+i}}\in\Q$ is a constant (see Proposition~\ref{prop:homog}).

From the Puiseux expansion we can deduce a so-called characteristic sequence. A classical result by Burau \cite{Bu} and Zariski \cite{Za} shows that
this sequence contains all the information about the topological type of the singularity. Mainly to fix the notation we recall the definition,
referring an interested reader to excellent books \cite{EN} or \cite{Wa}.

Assume that $q\ge p$. We define $q_0=0$, $p_0=p$. Assume inductively that $p_i$ and $q_i$ are chosen and $p_i>1$. We define
\[q_{i+1}=\min\{r\colon  c_r\neq 0,\, p_i\!\not|r\}\]
and $p_{i+1}=\gcd(p_i,q_{i+1})$. The procedure stops when $p_n=1$ (if it does never stop, then the parametrization \eqref{eq:param} is
not one to one). We define $(p;q_1,\dots,q_n)$ as the \emph{characteristic sequence} of the
singularity. We point out that $p_1,\dots,p_{n}$ are uniquely determined by $(p;q_1,\dots,q_n)$.

\begin{remark}
Our convention concerning the characteristic sequence is essentially that of \cite[Section 3]{Wa}.
\end{remark}

To fix the terminology, let us recall the following definition.

\begin{definition}
Let $c_jx^{j/p}$ be a term of the Puiseux expansion \eqref{eq:puis}. Let $m$ be the smallest index for which $p_m|j$. We shall call
$c_jx^{j/p}$ \emph{inessential} if either $m=0$ (so $p|j$) or there exists $j'<j$ such that $\gcd(p_{m-1},j')=p_m$ and $c_{j'}\neq 0$.
Otherwise it is called \emph{essential}. By a slight abuse of language we shall sometimes call $c_j$ or even $j$ itself essential or inessential.
\end{definition}

In other words, $c_j$ is inessential, if a change the value of $c_j$ has no effect on the topological type of the singularity. 
Conversely, if $c_j$ is an essential term, and we change the value of $c_j$ from
zero to non--zero, or from non--zero to zero, the topological type of the singularity is changed.

\begin{remark}
The notion 'essential' is not used in \cite{Wa}. In the Brieskorn and Kn\"orrer book 
only those terms are called essential, for which $c_j\neq 0$, see \cite[page 411]{BK}. 
\end{remark}

\begin{example}
For the singularity with $p=4$ and the Puiseux expansion $y=x^{3/2}+x^{9/4}$, the term $c_7x^{7/4}$ is essential even though $c_7=0$. Indeed, if we change
 $c_7$ to any non--zero value, then the characteristic sequence changes from $(4;6,9)$ to $(4;6,7)$. Similarly, the term $x^{3/2}$ is essential, because
if it is absent, then the characteristic sequence is $(4;9)$.
\end{example}

\subsection{Deformations}\label{ss:deform}

The main subject of the present article are deformations. 
We consider parametric deformations, i.e. we consider
the parametrization \eqref{eq:param} and assume that the coefficients $a_p,a_{p+1},\dots$ and $b_q,b_{q+1},\dots$ are continuous
functions of a parameter $s\in D$, where $D\subset\mathbb{C}$ is the unit disk. We assume that the functions $a_p(s)$ and $b_q(s)$
(the first coefficients of both equations in \eqref{eq:param}) do not vanish if $s\neq 0$.  Unless specified otherwise, a deformation
will always mean a parametric deformation in the sense we just specified. The value $s=0$ of the parameter will always be regarded
as special, corresponding to a 'most degenerate' singularity.

\begin{notation}\label{not:simple}
We will denote $a_j(s)$, $b_j(s)$, $c_j(s)$ and $\gamma_j(s)$ the corresponding
coefficients (or Puiseux coefficients) for a given value parameter $s$. However, when referring to $x$ or $y$, we shall use
the notation $x_s$ or $y_s$, in particular, for $s=0$, we shall write $x_0$, $y_0$, to avoid confusion with expressions like $x(t)$.

Furthermore, we consider the integers $p$, $q$ and $p'$ to be fixed. The integer $p$ is the order of $x_s$ at $t=0$ for $s\neq 0$, 
$p'$ denotes the order of $x_0$ and $q$ denotes $\ord_{t=0}y_s$ for $s\neq 0$.
\end{notation}

Parametric deformations are equinormalizable (see e.g. \cite[Section II.2.6]{GLS}, \cite{Za-book,Te}), because
one can regard the parametrization \eqref{eq:param} as a normalization of the singularity and therefore $\delta$--constant. 
One of the main, and the most difficult questions
in the theory is the adjacency problem. We can phrase it as follows.

\begin{question}
Given two singularities with characteristic sequences $\kappa=(p;q_1,\dots,q_n)$ and $\kappa'=(p';q_1',\dots,q_{n'}')$, does there exist a deformation
(in the sense described above) such that for $s\neq 0$ the singularity has characteristic sequence $\kappa$ and
for $s=0$ the singularity has characteristic sequence $\kappa'$ ?
\end{question}

\begin{remark}
In what follows, if such deformation exists, we shall say that $\kappa$ is \emph{adjacent} to $\kappa'$. 
\end{remark}

Although we restrict ourselves to a class of very explicit deformations, a complete answer to the above question seems to be extremely difficult 
and surprisingly
few full results (in the sense of conditions which are both sufficient and necessary) are known.
The only two completely solved cases are
$p'=p$ (see Section~\ref{s:completelysolved}) or $p=2$, $p'=3$, done by Petrov \cite{Pet} in a different language, see also \cite[Section 3]{BZ3}.
There are lots of different approaches known in the literature. 
The first one is the semicontinuity of the Milnor numbers (see e.g. \cite[Chapter 2]{Zol}). It is well known, see \cite{LR}, that if
under a deformation, the Milnor number is constant, then the deformation is topologically trivial. The semicontinuity of the Milnor numbers can
be generalized to the semicontinuity property of spectra \cite{Var}. Moreover, to obstruct adjacencies one can study dual graphs of
resolutions \cite{Ke}, Enriques diagrams \cite{ACR} or zero dimensional schemes (see \cite{GLS} and references therein). 
There are various approaches using knot theory, like \cite{Bo1} and \cite{Bo2}. 
One of the most promising is currently \cite{Ba}, relying on combinatorial properties of iterated torus knots.

The method we use in the paper is to look at the Puiseux polynomials (see Section~\ref{sec:puispoly}) from the point of view of
ordinary differential equations. 
In many cases we are able to restrict the possible Puiseux expansion in the limit. As we shall see, vanishing of \emph{inessential} terms when $s\neq 0$,
might imply vanishing of \emph{essential} terms in the limit $s=0$. 

We point out that studying properties of the Puiseux expansion under deformations is a method with a long history. Nevertheless, main applications concern
studying topologically equisingular deformations, from the point of view of analytic, not topological, 
equivalence see e.g. \cite{Za-book}. If we do not assume
topological equisingularity, the behaviour of the Puiseux expansion under deformation is much more complicated, compare Section~\ref{ss:example2to3}.

\subsection{Statement of results}
The first main result that we prove in this paper is the following.

\begin{theorem}\label{thm:first}
Suppose that we are given a deformation such that for $s\neq 0$ the singularity has the Puiseux expansion
\begin{equation}\label{eq:literally}
y_s=c_{r_0}(s)x_s^{r_0/p}+c_{r_1}(s)x_s^{r_1/s}+\dots,
\end{equation}
where
$r_1>r_0$ are positive integers. 
Suppose that for $s=0$ we have $\ord_{t=0}x_0=p'$. Then,
the first terms of the Puiseux expansion at $s=0$ are the following
\[y_0=c_{r'_0}x_0^{r_0'/p'}+c_{r_1'}x_0^{r_1'/p'},\]
where $r'_0=r_0p'/p$ and $r_1'\ge r_1+(p-p')$. Furthermore, if $r'_0\not\in\Z$, then $c_{r'_0}=0$.
\end{theorem}
\begin{remark}
Some authors, while writing a Puiseux expansion like \eqref{eq:literally} do not write inessential terms, even if they are not zero. 
This is not the case in this article. Equation~\eqref{eq:literally} should be read literally. The omitted terms $c_{r_0+1},\dots,c_{r_1-1}$ are all zero.
\end{remark}

The proof is given in Section~\ref{sec:prooffirst}.
In Section~\ref{sec:converse} we prove a partial converse to Theorem~\ref{thm:first}.

\begin{theorem}\label{thm:converse}
Let $p,r_0,r_1,p'$ be positive integers such that $r_1>r_0$, $p'>p$, $r_0<r_1+(p-p')$ and
\[r_1-r_0<\frac{p+r_0}{\gcd(p,r_0)}.\] 
Then, given any Puiseux expansion of the form
\begin{equation}\label{eq:referencePuis}
y=d_{r_1+p-p'}x^{(r_1+p-p')/p'}+d_{r_1+p-p'+1}x^{(r_1+p-p'+1)/p'}+\dots,
\end{equation}
and for an arbitrary $N\ge 0$, there exists a deformation such that for $s\neq 0$
the corresponding singularity has the expansion starting with
\[y_s=c_{r_0}x_s^{r_0/p}+c_{r_1}x_s^{r_1/p}+\dots,\]
where $c_{r_0}c_{r_1}\neq 0$, and
for $s=0$, the first $N$ terms of the Puiseux expansion agree with those of \eqref{eq:referencePuis}.
\end{theorem}

The methods used to prove Theorem~\ref{thm:first} can be generalized. The following result we shall prove in Section~\ref{sec:proofsecond}. 

\begin{theorem}\label{thm:second}
Let us consider a deformation, which for $s\neq 0$ has the Puiseux expansion
\[y_s=c_{r_0}(s)x_s^{r_0/p}+\dots+c_{r_n}(s)x_s^{r_n/p}+\dots,\]
where $r_n>r_{n-1}>\dots>r_0$ are positive integers.
Assume that the order of $x_0$ is $p'>p$. Suppose furthermore that for any $i=0,\dots,n-1$, either $r_i/p'\in\mathbb{Z}$, or $r_ip/p'\not\in\Z$. 
Then, up to a an overall change of coordinates of type  $(x,y)\to (x,y-Q(x))$, where $Q$ is a polynomial, we have
\[\ord y_0\ge r_n-(2n-1)(p'-p).\]
\end{theorem}
\begin{remark}
If $p$ and $p'$ are coprime, then for any $r_i\in\Z$, $r_i/p'\in\mathbb{Z}$ if and only if  $r_ip/p'\in\Z$.
\end{remark}

As we shall see in Section~\ref{sec:discussion}, Theorem~\ref{thm:second} is not optimal, i.e. the bounds for $y_0$ can 
sometimes be improved.

Apart from these theorems, we prove also two other results, related to the properties of functions $\gamma_{q+1},\gamma_{q+2},\dots$. 
The first one
is the non--genericity of the functions~$\gamma_{q+i}$, $i=1,2,\dots$ (see Theorem~\ref{thm:nongen} or Section~\ref{sec:bernkou}), which might be related to
the conjecture about semicontinuity of codimensions (see Section~\ref{s:csoc}). The other one is the relation between various derivatives
of the functions $\gamma_{q+1},\dots,$ encoded in WDVV equations, see Section~\ref{S:WDVV}.

\subsection{Conjectural semicontinuity of codimension}\label{s:csoc}
The content of this subsection serves as a motivation to methods, we shall be using in Section~\ref{sec:nongen}.

A codimension or $\ol{M}$-number was defined in \cite{Or} and independently in \cite{BZ1}. It was then studied in detail in \cite{BZ3}. For
cuspidal singularities it is defined as follows.

\begin{definition}
Let $\kappa=(p;q_1,\dots,q_n)$ be a characteristic sequence of the singularity. 
Let $p_1=p$ and for $j>1$ let $p_j=\gcd(p_{j-1},q_j)$. The \emph{codimension} of the singularity is defined as
\[\nu=p_1+q_1-\left\lfloor\frac{q_1}{p_1}\right\rfloor+\sum_{j=2}^{n-1}\left(q_j-q_{j-1}-\left\lfloor\frac{q_j-q_{j-1}}{p_j}\right\rfloor\right)-3.\]
\end{definition}

The motivation for this definition is the following. Consider the space $\C[t]\times\C[t]$ of parametrized curves in $\mathbb{C}^2$ (an
element $(x,y)\in\C[t]\times\C[t]$ defines a parametric curve in $\C^2$ by $t\to (x(t),y(t))$). For a given characteristic sequence $\kappa$,
let us consider $V_{\kappa}\subset\C[t]\times\C[t]$ consisting of those parametrized curves which for some $t_0$ (not necessarily $t_0=0$) has a cuspidal singularity
with the characteristic sequence $\kappa$. The codimension of $V_{\kappa}$ in $\C[t]\times\C[t]$ is exactly $\nu$ (see \cite{BZ3}). 

\begin{example}\label{ex:codim}
Consider $\kappa=(4;6,9)$. If $(x,y)\in V_{\kappa}$, then at some $t_0$, we can write $x(t)=x(t_0)+(t-t_0)^4a_4+\dots$, $y(t)=y(t_0)+(t-t_0)^4b_4+\dots$
(this gives $6$ conditions for vanishing of $a_1,a_2,a_3,b_1,b_2,b_3$). We substitute one because we are free to choose $t_0$. 
Moreover, either $a_4\neq 0$,
or $b_4\neq 0$. If $a_4\neq 0$, the Puiseux
expansion of $y-y(t_0)$ in the powers of $x-x(t_0)$ has form $y=c_4x+c_5x^{5/4}+c_6x^{6/4}+c_7x^{7/4}+c_8x^{8/4}+c_9x^{9/4}$ (if $a_4=0$ then $b_4\neq 0$
and we expand $x$ in powers of $y$ instead). We must have
$c_5=c_7=0$ and $c_6c_9\neq 0$. The last two conditions are open. The codimension is equal to $6-1+2=7$.
\end{example}

It is natural
to expect that if $\kappa$ and $\kappa'$ are two different characteristic sequences and $\nu(\kappa)\ge\nu(\kappa')$, then the closure of $V_\kappa$
in $\C[t]\times\C[t]$ omits $V_{\kappa'}$. Since $V_{\kappa'}$ has non--empty intersection with the closure of $V_\kappa$ if and only if $\kappa$
is adjacent to $\kappa'$, the above expectation is equivalent to the conjecture, that we now state.

\begin{conjecture}\label{conj:kappa}
If $\kappa$ is adjacent to $\kappa'$, then $\nu(\kappa)<\nu(\kappa')$.
\end{conjecture}

This conjecture was formulated independently by many authors, including Orevkov \cite{Or}, Christopher and Llynch \cite{CL} 
or Zoladek and the author \cite{BZ1,BZ3}. The validity of this conjecture would give very strong obstructions for the adjacency problem, and combined
with the current knowledge, it would give an answer close to the optimal (but see Example~\ref{ex:1016} and the discussion below it to see that Conjecture~\ref{conj:kappa}
is not enough itself).
The conjecture is still open and appears to be very difficult. 

Example~\ref{ex:codim} shows that $V_\kappa$ is given by equations of the type $\{a_i=0\}$, $\{b_j=0\}$ or $\{c_k=0\}$ (it is better to use $\gamma_k$ 
from \eqref{eq:cqi} instead
of $c_k$, because $c_k$ is in general a multivalued function), and some inequalities of type $\{a_i\neq 0\}$, $\{b_j\neq 0\}$ or $\{\gamma_k\neq 0\}$. 
If all $\gamma_k$'s were sufficiently generic polynomials, this would be
a major step towards proving Conjecture~\ref{conj:kappa}. However, in Section~\ref{sec:nongen} we prove a result which shows that the problem is
much more complicated. Informally, the result may be written as follows.

\begin{theorem}\label{thm:nongen}
The functions $\gamma_1,\dots,\gamma_k$ fail to satisfy suitably but naturally defined Bernstein--Kuschnirenko genericity conditions.
\end{theorem}

In Section~\ref{sec:bernkou} we explain in details, what Bernstein--Kuschnirenko genericity conditions do we have in mind. This does not mean
that Conjecture~\ref{conj:kappa} is false, but only that it is subtler than expected.

\section{Simple results about Puiseux expansions}
\subsection{Homogeneity of functions $\gamma_k$.}
We prove now a proposition, which we will use in Section~\ref{sec:nongen}. 
The result is almost obvious (but rather tedious), we give it for completeness of the exposition.

\begin{proposition}\label{prop:homog}
For any $i\ge 0$, the functions~$\gamma_{q+i}$ defined in \eqref{eq:cqi} 
are sums of monomials $\eta\cdot a_p^{k_p}\cdot\ldots\cdot a_{p+i}^{k_{p+i}}b_{q+j}$, where $\eta\in\mathbb{Q}$,
$k_p+\dots+k_{p+i}=i$, $j+k_{p+1}+2k_{p+2}+\dots+ik_{p+i}=i$.
\end{proposition}
\begin{proof}[Sketch of proof]
We proceed by induction on $l$. For $l=0,1,2$ the statement holds (see \eqref{eq:stupidgamma}). Assume that it holds for $l-1$.
Observe now that for $r\in\Z_{>0}$, the coefficient at $t^{r+l}$ in $a_p^{-r/p}x^{r/p}$ is a sum of terms of the form 
\begin{equation}\label{eq:form}
\eta\cdot \frac{a_{p+1}^{k_1}}{a_p^{k_1}}\cdot\ldots\cdot\frac{a_{p+l}^{k_l}}{a_p^{k_l}},
\end{equation}
where $\eta\in\Q$, $\sum jk_j=l$ and $\sum k_j\le l$. This follows immediately from the Taylor expansion 
\begin{multline*}
a_p^{-r/p}x^{r/p}=t^r(1+t\frac{a_{p+1}}{a_p}+t^2\frac{a_{p+2}}{a_p}+\dots)^{r/p}=\\
=t^r\left(1+\frac{r}{p}\left(t\frac{a_{p+1}}{a_p}+\dots\right)+\frac{r(r-p)}{2p^2}\left(t\frac{a_{p+1}}{a_p}+\dots\right)^2+\dots\right).
\end{multline*}
Now the inductive assumptions together with the above observations can be used to show that the coefficient at $t^{q+l}$ of
\[y-(c_qx^{q/p}+\dots+c_{q+l-1}x^{(q+l-1)/p})\]
is of the form 
\begin{equation}\label{eq:tobecancelledbycqi}
b_{q+l}-\sum \eta_{k_p,\dots,k_{p+l}}b_{q+l-j}\frac{a_{p+1}^{k_1}}{a_p^{k_1}}\cdot
\ldots\cdot\frac{a_{p+l}^{k_{l}}}{a_p^{k_{l}}},
\end{equation}
where we sum over non--negative integers $k_p,\dots,k_{p+l}$ such that $0\le \sum k_i\le l$ and $j=\sum ik_{p+i}$ is in the range $1,\dots,l$.
Here $\eta_{k_p,\dots,k_{p+l}}$ is a rational number and $\eta_{0,\dots,0,1}=\frac{q}{p}$.
But the expression \eqref{eq:tobecancelledbycqi} is also equal to $c_{q+i}a_p^{(q+i)/p}$. The induction step follows.
\end{proof}

\subsection{Boundedness of $c_{q+i}(s)$ and its consequences.}\label{ss:example2to3}
As we already mentioned in the introduction,
when studying deformations we simply allow the coefficients $a_p,\dots$ and $b_q,\dots$ of \eqref{eq:param} to vary with a parameter $s$.
Then $c_{q+i}$ and $\gamma_{q+i}$ can be regarded as functions of $s$.

We begin with a specific example of a deformation, when we show that for certain values of $i$, the function $c_{q+i}(s)$ 
diverges to infinity as $s\to 0$, and explain what happens
if we can ensure boundedness. 

Let us consider the deformation given by
\begin{equation}\label{eq:def10I}
\begin{split}
x_s(t)&=s t^2+t^3\\
y_s(t)&=b_2(s)t^2+b_3(s)t^3+\dots,
\end{split}
\end{equation}
Assume that $b_i$'s are chosen in such a way that for $s\neq 0$ the resulting singularity is an
$A_{10}$ singularity (the characteristic sequence is $(2;11)$). This amounts to the fact that we have
\begin{equation}\label{eq:def10II}
y_s=c_2(s)x_s+c_4(s)x_s^2+
c_6(s)x_s^3+c_8(s)x_s^4+c_{10}(s)x_s^5+c_{11}(s)x_s^{11/2}+\dots.
\end{equation}
Substituting $x$ from \eqref{eq:def10I} into \eqref{eq:def10II} yields
\begin{equation}\label{eq:def10III}
\begin{split}
y_s=&t^2 sc_2(s)+t^3c_2(s)+t^4s^2c_4(s)+t^5\cdot 2sc_4(s)+\\
&+t^6(c_4(s)+c_6(s)s^3)+t^7\cdot 3s^2c_6(s)+\\
&+t^8(3sc_6(s)+s^4c_8(s))+t^9(c_6(s)+4s^3c_8(s))+\\
&+t^{10}(6sc_8(s)+s^5c_{10}(s))+\dots
\end{split}
\end{equation}

Then we have the following observation.

\begin{lemma}\label{lem:b11}
If all $c_2,c_4,c_6,c_8$ and $c_{10}$ remain bounded from above while $s\to 0$ then the singularity at $s=0$ has characteristic sequence $(3;b)$
with $b\ge 11$.
\end{lemma}
\begin{proof}
By assumptions, when $s\to 0$, all terms in \eqref{eq:def10III} with $s$ in positive power converge to $0$. Hence $y_0$ has order at least $11$.
\end{proof}

Example~\ref{ex:specb8} shows that $(2;11)$ is adjacent to $(3;8)$ hence the assumptions of Lemma~\ref{lem:b11} are not always true. 
What we can prove is the following.
\begin{lemma}\label{l:2.5}
$sc_6(s)$ is bounded as $s\to 0$.
\end{lemma}
\begin{proof}
The relation between $c_6(s),c_8(s),b_8(s)$ and $b_9(s)$ can be written as
\[
\begin{pmatrix}
3 &1\\
1&4
\end{pmatrix}
\left(
\begin{matrix}
s c_6(s)\\ s^4 c_8(s)
\end{matrix}
\right)=
\begin{pmatrix}
b_8(s)\\ sb_9(s)
\end{pmatrix}.
\]
The vector on the right hand side is bounded as $s\to 0$. As the determinant of the matrix is non--zero, we infer that $sc_6$ and $s^4c_8$ are bounded.
\end{proof}

From Lemma~\ref{l:2.5} it follows that in the limit expansion $b_7(s)\to 0$. This implies the following result.
\begin{corollary}\label{cor:specb8} 
If $(2;11)$ is adjacent to $(3;b)$, then $b\ge 8$.
\end{corollary}

We emphasize that non--degeneracy of the matrix $\left(\begin{smallmatrix} 3 & 1\\ 1& 4\end{smallmatrix}\right)$ is a key point in the proof of 
Corollary~\ref{cor:specb8}. 
More complicated matrices appear when one tries to study possible adjacencies of $(p;q)$ to $(p+1;q')$ for general $p,q,q'>1$. 
Finding a general formula for the determinants (to show that they are non-zero)
is difficult, but doable, at least in theory. We plan to investigate such adjacencies in a future paper.

\begin{example}\label{ex:specb8}
The bound $b\ge 8$ is optimal. Indeed, let us consider the 
family $x_s(t)=st^2+t^3$, $y_s(t)=4s^2t^6+12st^7+t^8$. For $s\neq 0$ the singularity at the origin is an $A_{10}$ singularity,
while for $s=0$ we have $x(t)=t^3$, $y(t)=t^8$.
\end{example}

A few important conclusions are in order. 
\begin{itemize}
\item[(1)] The Puiseux coefficients $c_{q+i}(s)$ might diverge as $s\to 0$. 
\item[(2)] The polynomials $\gamma_{q+i}(s)$ are obviously convergent, but if $a_p(0)=0$, the meaning of $\gamma_{q+i}(0)$ is
unclear. Definitely it is not the polynomial $\gamma_{q+i}$ of the singularity at $s=0$. In fact, the order of $x$ (equal to $p$ for $s\neq 0$)
enters the definition of $\gamma_{q+i}$ (see \eqref{eq:stupidgamma} for instance), but if $a_p(0)=0$, the order of $x_0$ (cf. Notation~\ref{not:simple})
is bigger than $p$.
\item[(3)] Boundedness, or vanishing of inessential  Puiseux coefficients when $s\neq 0$ might influence the limit (in the
above case $c_2,c_4,\dots,c_{10}$ are inessential if $s\neq 0$).
\item[(4)] In the limit, some inessential terms might vanish as well. For example, in Lemma~\ref{lem:b11} we proved that $\ord y_0\ge 11$. This is
more than just the statement that $x_0,y_0$ has singularity of type $(3;b)$, we know that in the expansion 
$y_0=c_3x_0^{3/3}+c_6x_0^{6/3}+c_9x_0^{9/3}+c_{11}x_0^{11/3}+\dots$, the terms $c_3$, $c_6$ and $c_9$ are equal to $0$.
\end{itemize}

In other words, the inessential Puiseux terms might play a role in determining the type of the limit singularity. We will see more examples in
Section~\ref{sec:ODE}.

\subsection{Adjacency in the case $p=p'$.}\label{s:completelysolved}
The case, when $a_p(0)\neq 0$ can be solved completely. To describe the solution, we introduce some additional notation:

For given $p>0$ we consider the space $\C_p[t]\subset\C[t]$ of polynomials with terms up to $t^{p-1}$ inclusively equal to $0$. For a given
singularity with characteristic sequence $\kappa=(p,q_1,\dots,q_n)$, we consider the set $V^p_{\kappa}\subset\C_p[t]\times\C_p[t]$ of pairs
of polynomials $(x,y)$ such that the curve $t\to(x(t),y(t))$ has singularity at $t=0$ with characteristic sequence $\kappa$, and $\ord x=p$
(i.e. $a_p\neq 0$). The set $V^p_{\kappa}$ is given by a set of equations of type $c_{j}=0$ and inequalities $c_{k}=0$.

\begin{definition}\label{def:I}
The sets $I_{\kappa}$ and $J_{\kappa}$ are defined as finite sets such that
\[V^p_{\kappa}=\bigcap_{i\in I_{\kappa}}\{c_i=0\}\cap\bigcap_{j\in J_{\kappa}}\{c_j\neq 0\}.\]
\end{definition}
$I_\kappa$ is the set of indices of vanishing essential Puiseux coefficients, $J_\kappa$ is the set of non--vanishing essential Puiseux terms.
In particular, if $\kappa=(p;q_1,\dots,q_n)$, then $J=(q_1,\dots,q_n)$ and the largest element in $I_\kappa$
is smaller than $q_n.$

The definition being slightly artificial, we give some simple examples.

\begin{example}
If $\kappa=(4;6,9)$, then $I_{\kappa}=(5,7)$ and $J_{\kappa}=(6,9)$. 
If $\kappa=(7;11)$, then $I_{\kappa}=(8,9,10)$. For $\kappa=(2;2k+1)$ we have $I_{\kappa}=(3,5,\dots,2k-1)$.
\end{example}

We remark that $\#I_{\kappa}+2p-3=\nu(\kappa)$ is the codimension. We have a following lemma, which we will not need until Section~\ref{sec:whatisl0}.

\begin{lemma}\label{lem:semigroup}
The set $\Gamma_\kappa=\{n\colon p+n\not\in I_{\kappa}\}$ is a semigroup.
\end{lemma}

\begin{remark}
$\Gamma_\kappa$ \emph{is not} the semigroup of the singularity. For example, 
for $\kappa=(4;6,7)$ we have $I_{\kappa}=(5)$, $\Gamma_{\kappa}=(2,3,\dots)$, but the semigroup
of the singularity is generated by $(4,6,13)$ by e.g. \cite[Theorem 4.3.5]{Wa}.
\end{remark}

\begin{proof}[Proof of Lemma~\ref{lem:semigroup}]
Let $\kappa=(p;q_1,\dots,q_n)$ and  $p_0,p_1,\dots,p_n$ be as Section~\ref{sec:prel}.
Let us take two integers $a,b\ge 0$ such that $a+p,b+p\not\in I_\kappa$. We need to show that $a+b+p\not\in I_\kappa$, either.

\emph{Claim.} Given an integer $a>0$, if $m$ is the smallest possible index such that $p_m|a$ and $a+p\not\in I_{\kappa}$, then $a+p\ge q_m$.

To prove the claim observe, that if $a+p\not\in I_{\kappa}$ then $a+p$ is inessential or non--vanishing. 
But this means that it must be a non--vanishing essential term $c_{a'+p}$ in the Puiseux expansion before $c_{a+p}$ (i.e. $a'\le a$)
for which $p_m|a'$. Let us take the smallest such $a'$. But then, by the very definition, $p+a'=q_m$.

Given the claim let us choose smallest indices $m,m'$ such that $p_m|a$ and $p_{m'}|b$. Assume that $m\ge m'$, 
then $a+p\ge q_m$, $b+p\ge q_{m'}$ so $a+b+p>q_{m}$ and $p_m|(a+b+p)$,
hence $c_{a+b+p}$ is inessential (because all terms after $q_m$, which are divisible by $p_m$ are essential). In particular $a+b+p\not\in I_\kappa$.
\end{proof}

Now we are ready to state the main result of this subsection. We expect it to be standard and well-known to the experts, but we
could not find a precise reference.

\begin{proposition}
Let $\kappa=(p;q_1,\dots,q_n)$ and $\kappa'=(p';q_1',\dots,q_{n'}')$. If $p=p'$ then $\kappa$ is adjacent to $\kappa'$ if and only if
$I_{\kappa}\subset I_{\kappa'}$.
\end{proposition}
\begin{proof}
\underline{The 'only if' part.}
Let us choose a change of variables near $(t,s)=(0,0)$ given by $(t,s)\to (\tau,s)=(x_s^{1/p},s)$. Since $a_p(0)\neq 0$, this
is a change of variables analytic in $t$ and continuous in $s$. In the new variable we have
\begin{align*}
x_s(\tau)&=\tau^p\\
y_s(\tau)&=c_p(s)\tau^q+c_{p+1}(s)\tau^{q+1}+\dots.
\end{align*}
The functions $c_{k}(s)$ are the same as above, since they do not depend on the parametrization of the curve. In particular, they
are continuous. Now if $i\in I_{\kappa}$, then for $s\neq 0$ we have $c_i(s)=0$. By continuity $c_i(0)=0$. Assume there exists 
$i_0\in I_{\kappa}\setminus I_{\kappa'}$ and let us take the smallest such $i_0$. Clearly, $i_0\not\in J_{\kappa'}$, because $c_{i_0}(0)=0$.
Let $\kappa=(p,q_1,\dots,q_n)$ and $\kappa'=(p,q_1',\dots,q'_{n'})$ and let $p_1,\dots,p_n$ (respectively $p_1',\dots,p'_{n'}$)
be as in Section~\ref{sec:prel}. Let $m$ and $m'$ be such that $q_m<i_0<q_{m+1}$ and $q_{m'}'<i_0<q'_{m'+1}$
(if $i_0>q_n$ we put $m=n$, if $i_0>q'_{n'}$ we define $m'=n'$). 

Since $i_0\in I_\kappa$, we infer that $p_m$ does not divide $i_0$, otherwise $i_0$ is inessential for $s\neq 0$ and does not belong to $I_\kappa$.
As for any $j$, the condition $c_j(s)=0$ for $s\neq 0$ implies $c_j(0)=0$, $p_{m}$ divides $p_{m'}'$. Indeed,  by construction 
$p_m$ is the greatest common divisor of all $j<q_{m+1}$ such that $c_j(s)\neq 0$ for $s\neq 0$. Therefore $p_{m'}'$ does
not divide $i_0$. As $i_0$ is not equal $q'_j$ for any $j$, we infer that $i_0$ is essential for $s=0$, hence $i_0\in I'$, contradiction. 

\smallskip
\underline{The 'if' part.} We shall present an explicit deformation of $\kappa$ to $\kappa'$. Let us define a deformation
\begin{align*}
x_s(t)&=t^p\\
y_s(t)&=\sum_{j=p}^\infty b_j(s)t^j,
\end{align*}
where
\[
b_j(s)=
\begin{cases}
0& \textrm{ if $j\in I$}\\
s& \textrm{ if $j\in I'$ and $j\not\in I$}\\
1& \textrm{ otherwise.}
\end{cases}
\]
For $s\neq 0$ the singularity is $\kappa$, for $s=0$, the singularity is $\kappa'$.
\end{proof}

\section{Ordinary differential equations related to ODE's}\label{sec:ODE}

\subsection{Puiseux polynomials}\label{sec:puispoly}
For further applications, we review the construction from \cite{Bo3}. Let us be given a singularity as in \eqref{eq:param}. We shall write
the Puiseux expansion \eqref{eq:puis} in the following form.
\begin{equation}\label{eq:puis2}
y=c_{r_0}x^{r_0/p}+c_{r_1}x^{r_1/p}+c_{r_2}x^{r_2/p}+\dots.
\end{equation}
Here $r_0=q$ and $r_1<r_2<\dots$ are chosen so that $c_{r_j}\neq 0$ (we begin the notation with $r_0$, not with $r_1$ to agree with
\cite{Bo3}). The sequence $(r_0,\dots,r_n,\dots)$ together with the value $p$ determines the
topological type of the singular point, but it contains more information: we know which \emph{inessential} Puiseux coefficients vanish. Surprisingly,
the discussion below depends also on this additional data.

In what follows $\dot x$ and $\dot y$ denote derivatives over $t$, $\ddot x$ is the second derivative. Let us define
\begin{equation}\label{eq:p1}
P_1=\dot yx-\frac{r_0}{p}\dot x y,
\end{equation}

Then we have the following result
\begin{lemma}[see \expandafter{\cite[Lemma 1]{Bo3}}]\label{lem:ordp1}
The order of $P_1$ at $t=0$ is equal to $r_1+(p-1)$.
\end{lemma}

We can define $P_2,\dots,P_n,\dots$ inductively by
\begin{equation}\label{eq:pk}
P_{k+1}=x\dot x\frac{d}{dt}P_k-\left(\frac{r_k}{p}{\dot x}^2+(2k-1)\ddot x x\right)P_k.
\end{equation}
We have
\begin{lemma}[see \expandafter{\cite[Lemma 2]{Bo3}}]\label{lem:ordp2}
The order of $P_k$ at $t=0$ is equal to $r_k+(2k-1)(p-1)$.
\end{lemma}

The functions $P_1,\dots,$ are called the \emph{Puiseux polynomials}. We shall show that they behave quite well under deformations.

\subsection{Proof of Theorem~\ref{thm:first}.}\label{sec:prooffirst}

For the convenience of the reader we recall the statement of Theorem~\ref{thm:first}.

\begin{theorem}
Suppose that  we are given a deformation such that for $s\neq 0$ the singularity has the Puiseux expansion
\[y_s=c_{r_0}x_s^{r_0/p}+c_{r_1}x_s^{r_1/s}+\dots,\]
where $r_1>r_0$ are positive integers. Suppose that at $s=0$ we have $\ord_{t=0}x_0=p'$. Then,
the first terms of the Puiseux expansion at $s=0$ are the following
\[y_0=c_{r'_0}x_0^{r_0'/p'}+c_{r_1'}x_0^{r_1'/p'},\]
where $r'_0=r_0p'/p $ and $r_1'\ge r_1+(p-p')$. Furthermore, if $r'_0\not\in\Z$, then $c_{r'_0}=0$.
\end{theorem}
\begin{proof}
For the given deformation, we write $P_1(s),P_2(s)$
for the corresponding Puiseux polynomials. To avoid ambiguities we shall also write
\begin{equation}\label{eq:rkpk0}
S_k:=P_k(0),
\end{equation}
which is a function only of the variable $t$.
As polynomial functions of $a_p(s),\dots$ and $b_q(s),\dots$,
the coefficients of $P_1(s),\dots,P_k(s),\dots$ are also continuous in $s$. In particular for any $k=1,\dots$ we have

\begin{equation}\label{eq:ordsemic}
\ord_{t=0}S_k\ge \lim\sup_{s\to 0}\ord_{t=0}P_k(s).
\end{equation}
\begin{remark}
In general $S_k$ is not the Puiseux polynomial for the singularity at $s=0$, compare with item (2) in the list after Example~\ref{ex:specb8}.
\end{remark}

By assumption of the theorem and Lemma~\ref{lem:ordp1} we have for $s\neq 0$
\[\ord_{t=0} P_1(s)=r_1+(p-1).\]
Then 
\begin{equation}\label{eq:neweq}
\ord_{t=0} S_1\ge r_1+(p-1).
\end{equation}
Observe that $S_1$ still satisfies 
\begin{equation}\label{eq:ODEforS}
\dot y_0x_0-\frac{r_0}{p}y_0\dot x_0=S_1, 
\end{equation}
which we now regard --- and this is the key point of the
proof --- as a differential equation with known functions $x_0$ and $S_1$ and unknown $y_0$. Solving it we obtain
\begin{equation}\label{eq:solvey0}
y_0=Cx_0^{r_0/p}+x_0^{r_0/p}\int^t \frac{S_1}{x_0^{r_0/p+1}},
\end{equation}
where $C\in\C$ is the integration constant.
(see Remark~\ref{rem:integral} below for explanation of the symbol $\int^t$).
The expression
\[x_0^{r_0/p}\int^t\frac{S_1}{x_0^{r_0/p+1}}\]
has order equal to $\ord_{t=0}S_1-(p'-1)$, hence it can be expanded in fractional powers of $x_0$ as
\[c_{r_1'}x_0^{r_1'/p'}+c_{r_1'+1}x_0^{(r_1'+1)/p}+\dots,\]
where $r_1'=\ord_{t=0}S_1-p'-1\ge r_1+p-p'$. On substituting this into \eqref{eq:solvey0} we conclude the first part of the proof.

Of $r_0'\not\in\Z$, then $x_0^{r_0'/p'}$ is multivalued near $t=0$. If $C\neq 0$, then the term $Cx_0^{r_0/p}=Cx_0^{r_0'/p'}$
is the only multivalued term of both sides of \eqref{eq:solvey0}, which is absurd. Hence, if $r_0'\not\in\Z$, we have $C=0$.
\end{proof}
\begin{remark}\label{rem:integral}
If $\zeta(t)$ is given by a series $\sum_{n\in\Z}a_nt^{n/p}$ (with $a_{-p}=0$), then the integral $\int^t\zeta(t)$ is defined as 
$\sum_{n\in\Z}\frac{p}{n+p}a_nt^{1+n/p}$. If the integrand is an analytic multivalued function near $t=0$, 
then the integral is also analytic and multivalued. It is easy to check, that with this definition $y_0$ indeed satisfies \eqref{eq:ODEforS}.
\end{remark}

\subsection{Proof of Theorem~\ref{thm:second}}\label{sec:proofsecond}
Let us recall the statement of Theorem~\ref{thm:second}.

\begin{theorem}
Let us consider a deformation, which for $s\neq 0$ has the Puiseux expansion
\[y_s=c_{r_0}x_s^{r_0/p}+\dots+c_{r_n}x_s^{r_n/s}+\dots,\]
where $r_n>r_{n-1}>\dots>r_0$ are positive integers.
Assume that the order of $x_0$ is $p'>p$. Suppose furthermore that for any $i=0,\dots,n-1$, either $r_i/p'\in\mathbb{Z}$, or $r_ip/p'\not\in\Z$. 
Then, up to an overall change of coordinates of type  $(x,y)\to (x,y-Q(x))$, where $Q$ is a polynomial, we have
\[\ord y_0\ge r_n-(2n-1)(p'-p).\]
\end{theorem}

\begin{proof}
Let us consider the $n-$th Puiseux polynomial $P_n$. By Lemma~\ref{lem:ordp2}, for $s\neq 0$, we have
\[\ord_{t=0} P_n(s)=r_n+(2n-1)(p-1).\]
By passing to the limit $s\to 0$ we obtain
\[\ord_{t=0} S_n\ge r_n+(2n-1)(p-1),\]
where $S_n$ is defined in \eqref{eq:rkpk0}.
By induction we want to show that for $1\le k\le n$ we have
\begin{equation}\label{eq:degpk0}
\ord_{t=0} S_k\ge r_n+(2n-1)(p-1)-(2n-2k)(p'-1).
\end{equation}
Assume $1\le k\le n$.
The induction step is done by considering \eqref{eq:pk} as a non--homogeneous linear ODE with unknown $S_k$.
When $s=0$, we can rewrite \eqref{eq:pk} as:
\[\frac{S_k'}{S_k}=\frac{r_k}{p}\frac{\dot x_0}{x_0}+(2k-1)\frac{\ddot x_0}{\dot x_0}+\frac{S_{k+1}}{x_0\dot x_0}.\]
The general solution has the form
\[S_k=x^{r_k/p}\dot x^{2k-1}\left(C+\int^t\frac{S_{k+1}}{x^{1+r_k/p}\dot x^{2k}}\right),\]
where $C$ is the integration constant.
If $p'r_k/p\not\in\Z$, then  $x^{r_k/p}\dot x^{2k-1}$ is not analytic near $t=0$. Hence $C=0$,
otherwise $S_k$ is not analytic. 

If $p'r_k/p\in\Z$, then by assumptions we have $r_k/p\in\Z$. By \cite[Equations (6) and (7)]{Bo3} (it is also a straightforward consequence of
\eqref{eq:pk}) a change of variables of type $y\to y-C_1x^{r_k/p}$ induces, for any $m\ge 2$ a change $S_m\to S_m-\delta_mC_1x^{r_k/p}\dot x^{2m-1}$,
where $\delta_m=\prod_{j=1}^{m-1}\frac{r_k-r_j}{p}$. In particular $\delta_k\neq 0$. Hence a change of variables $y\to y-\frac{C}{\delta_k} x^{r_k/p}$
kills the term $C x^{r_k/p}\dot x^{2k-1}$ in $S_k$ and we can assume again that $C=0$.

Comparing now the orders of both sides, we finish the induction step.
In particular we know that
\[\ord S_1\ge r_n+(2n-1)(p-p')+(p'-1).\]
Arguing as in the proof of Theorem~\ref{thm:first} (see Section~\ref{sec:prooffirst}) we conclude that
$\ord y_0\ge r_n+(2n-1)(p-p')$.
\end{proof}

\subsection{Discussion of Theorems~\ref{thm:first} and \ref{thm:second}}\label{sec:discussion}

First we shall show an application of Theorem~\ref{thm:first}.

\begin{example}\label{ex:1016}
Suppose that we are given a deformation such that for $s\neq 0$, the singularity has Puiseux expansion
\[y_s=c_9(s)x_s+c_{17}(s)x^{17/9}+\dots\]
with $c_{17}(s)\neq 0$,
i.e. it is a $(9;17)$ singularity. If $\ord_{t=0} x_0=10$, then $\ord_{t=0}y_0\ge 16$. In other words, $(9;17)$ can be adjacent to $(10;16,17)$,
but is not adjacent to $(10;15,n)$. for any $n\ge 16$.
\end{example}

Theorem~\ref{thm:converse} allows us to construct an explicit deformation of the singularity $(9;17)$ to $(10;16,17)$, see Example~\ref{ex:1016again}.

\begin{remark}
The Milnor number of $(9;17)$ is $128$, whereas of $(10;15,16)$ is $130$. The semicontinuity of Milnor numbers does not obstruct the adjacency
of $(9;17)$ to $(10;15,n)$.
\end{remark}
\begin{remark}
The codimension of $(10;15,n)$ can be arbitrary large (it grows like $\frac45n$), so the conjectured semicontinuity of codimensions
(see Conjecture~\ref{conj:kappa}) is not
sufficient to obstruct the adjacency of $(9;17)$ to $(10;15,n)$ for large $n$.
\end{remark}
\begin{remark}
The Levine--Tristram signatures can be used to  obstruct the adjacency of $(9;17)$ to $(10;15,16)$  \cite{Bo1,Bo2}. Indeed, the value of Levine--Tristram
signature at $z=e^{2\pi i\cdot 0.165}$ for the torus knot $T_{9,17}$ is less than the corresponding value for the iterated torus knot 
$T_{5,8;81,2}$ (i.e. a $(81,2)$
cable on $T_{5,8}$, which is the link of the singularity $(10;15,16)$, see \cite{EN}). 
Adjacency would imply that the value is not less (this follows from \cite[Proposition 4.3 and Proof of Theorem 5.2]{Bo2},
altough this result is not written explicitely). 
However this obstruction does not prohibit the adjacency of $(9;17)$ to $(10;15,n)$ for $n\ge 17$.
\end{remark}

The next example shows the usage of Theorem~\ref{thm:second}. In general, the estimates are not so good, as one would expect.

\begin{example}
Let $N>0$ and consider a deformation with $p=16$ and the Puiseux expansion of a generic fiber
\[y_s=c_{16}(s)x_s^{16/16}+c_{20}(s)x_s^{20/16}+\dots+c_{4N+16}(s)x_s^{1+N/4}+c_{4N+17}(s)x_s^{(4N+17)/16}+\dots.\]
Suppose that the order of $x_0$ is $17$. Then by Theorem~\ref{thm:second}
\[\ord y_0\ge 17+2N.\]
It is interesting to compare the Milnor numbers. The general member $s\neq 0$ has
 $\mu_s=282+12N$, while $(17,17+2N)$ has $256+16N$. For large $N$, Theorem~\ref{thm:second} 
is a stronger obstruction than the semicontinuity of  Milnor numbers. However, for small $N$ it is much weaker.
\end{example}

When applying Theorem~\ref{thm:second} one should remember that using large values of  $n$ is not necessarily optimal, in fact, if $r_{n+1}-r_n<p'-p$, increasing
$n$ leads to \emph{weaker estimates}.

\begin{example}
Consider a deformation with Puiseux expansion
\[y_s=c_{16}x_s^{16/16}+c_{24}x_s^{24/16}+c_{32}x_s^{32/16}+c_{40}x_s^{40/16}+c_{42}x_s^{42/16}+c_{43}x_s^{43/16}\dots,\]
i.e. $r_0=16$, $r_1=24$, $r_2=32$, \dots. The Puiseux coefficient may depend on $s$. 
Assume that $\ord_{t=0}x_0=19$, so $p'-p=3$. Applying Theorem~\ref{thm:second} for $n=2$,
yields $\ord_{t=0}y_0\ge 23$, for $n=3$ we have $\ord_{t=0}y_0\ge 25$, for $n=4$ however we obtain $\ord_{t=0}y_0\ge 21$. For $n=5$ we get 
$\ord_{t=0}y_0\ge 16$.
\end{example}

The above example shows that there are lots of subtleties related to orders of polynomials $S_k$. The fact that looking at orders
of $S_4$ yields a weaker estimate than looking at orders of $S_3$ indicates, that
\[\ord_{t=0}S_4>\lim\sup_{s\to 0}\ord_{t=0}P_4(s).\]
Therefore there are some jumps of orders, which we still do not understand. A possible way to deal with that problem is to apply tools
from qualitative theory of linear non--homogeneous ODE's of Fuchs type (like dependence of solutions on a parameter $s$), but most
of the theory of such equations deals with homogeneous ODE's (see e.g. \cite[Chapter 8]{Zol}).

\section{Subsets defined by zeros of Puiseux coefficients}\label{sec:nongen}

\subsection{Setup}
Up to now, we were considering deformations as families of parametrization $t\to (x_s(t),y_s(t))$ depending on a deformation parameter $s$. Now we shall
slightly change a point of view. A deformation will be regarded as a (germ of) a complex curve in a space of coefficients of $(x(t),y(t))$ 
(see \eqref{eq:param}). We need first to define the space of coefficients. Taking a space of all parametrization \eqref{eq:param} is not the best
choice, the space would be infinite dimensional. We shall use the fact that each isolated singularity depends only on a finite Taylor expansion
(this follows from e.g. the Tougeron lemma, see \cite[Section 2.1]{Zol}) and restricts ourselves to polynomials of sufficiently high degree.

Let us fix $p$ and $q$ and choose a large enough integer $M$. The spaces $V_x$ and $V_y$ are defined as the linear spaces of coefficients respectively
$(a_p,a_{p+1},\dots,a_{p+M})$ and $(b_q,b_{q+1},\dots,b_{q+M})$. A point $(x,y)=(a_p,\dots,a_{p+M};b_q,\dots,b_{q+M})\in V_x\oplus V_y$ will be regarded as a 
cuspidal singularity parametrized locally by

\begin{equation}\label{eq:elementofV}
\begin{split}
x(t)&=a_pt^p+\dots+a_{p+M}t^{p+M}\\
y(t)&=b_qt^q+\dots+b_{q+M}t^{q+M}.
\end{split}
\end{equation}

Furthermore we shall introduce the space $\PV=\PP V_x\oplus\PP V_y\cong \C P^{M-1}\times\C P^{M-1}$ of projectivised parameters 
The main advantage
is that $\PV$ is compact. The homogeneity property (Proposition~\ref{prop:homog}) shows that
the functions $\gamma_{q+1},\gamma_{q+2},\dots$ give rise to homogeneous functions on $\PV$.
In particular we can look at their zero sets.

\begin{remark}\label{rem:elementofV}
An element $(x,y)\in\PP V$ shall be regarded as a parametrization \eqref{eq:elementofV} defined up to a multiplication of $x$ by a constant from $\C^*$
and $y$ by another constant.
\end{remark}

Consider a singularity $\kappa=(p;q_1,\dots,q_n)$ where $q_1=q$ and $q_n<q+M$. Let $I_{\kappa}$ be as defined in Definition~\ref{def:I}.
We define
\[I_\kappa^q=I_\kappa\cap[q,\infty).\]
Let
\[\tR_\kappa=\bigcap_{i\in I^q_\kappa}\{\gamma_i=0\}.\]
We would like to call $\tR_\kappa$ the defining set of the singularity.
We have a disappointing fact.

\begin{lemma}
For any $i\ge 0$, the set $\{a_p=0\cap b_qa_{q+1}=0\}\subset\PP V$ is a subset of $\{\gamma_{q+i}=0\}$.
\end{lemma}
\begin{proof}
We will show that $\gamma_{q+i}=0$ whenever $a_p=0$ and $a_{p+1}b_q=0$. By Proposition~\ref{prop:homog} all
monomials entering in $\gamma_{q+i}$ are of the form
\[b_{q+j}a_p^{k_p}a_{p+1}^{k_{p+1}}\dots a_{p+i}^{k_{p+i}},\]
where 
\begin{align*}
k_p+\dots+k_{p+i}&=i\\
j+k_{p+1}+2k_{p+2}+\dots+ik_{p+i}&=i.
\end{align*} 
If $k_p>0$, then the monomial vanishes at $\{a_p=0\}$.
So let us suppose that $k_p=0$.
Taking the difference of the two equations we arrive at
\[j+k_{p+2}+2k_{p+3}+\dots+(i-1)k_{p+i}=0.\]
As all the entries are non--negative, we infer that $j=k_{p+2}=\dots=k_{p+i}=0$. This corresponds to the monomial $b_qa_{p+1}^{i}$, which
vanishes on the set $\{b_qa_{p+1}=0\}$.
\end{proof}

It follows that if $\# I^q_\kappa\ge 3$, then $\tR_\kappa$ has not even correct codimension (i.e. equal to $\# I^q_\kappa$)! Furthermore, as $\gamma_{q+j}$ restricted
to $\{a_p=0\}$ is proportional to $b_qa_{p+1}^j$, the scheme structure of $\tR_\kappa$ near $\{a_p=0\}\cap\{b_qa_{p+1}=0\}$ is, in general, very
complicated. To remedy to this problem, we shall introduce the following definition, which we shall use throughout the section.

\begin{definition}\label{def:LI}
Let $I$ be a finite set of numbers from $\{q+1,\dots,q+M\}$, the scheme $L_I\subset\PP V$ is defined as
\[ L_I=\ol{\bigcap_{i\in I}\{\gamma_i=0\}\setminus\left\{a_p=0\right\}}.\]
If $I=I^q_{\kappa}$ for some singularity $\kappa$, then we shall write $L_\kappa$ instead of $L_{I^q_\kappa}$.
\end{definition}

The study of properties of $L_I$ is rather difficult, but not so difficult as $\tR_\kappa$. The following simple result is very helpful.
\begin{lemma}\label{lem:liissmooth}
For any $I$, the set $L_I$ (or $\tR_I$) is smooth away from the hypersurface $\{a_p=0\}$.
\end{lemma}
\begin{proof}
The proof is straightforward, even if not very enlightening.
Let $I=(i_1,\dots,i_\nu)$, where $i_1,\dots,i_\nu>0$. Consider the matrix of partial derivatives
\[D\gamma=\begin{pmatrix}
\frac{\partial{\gamma_{i_1}}}{\partial a_p}&\dots&\frac{\partial{\gamma_{i_1}}}{\partial a_{p+M}}&
\frac{\partial{\gamma_{i_1}}}{\partial b_q}&\dots&\frac{\partial{\gamma_{i_1}}}{\partial b_{q+M}}\\
\dots\\
\frac{\partial{\gamma_{i_\nu}}}{\partial a_p}&\dots&\frac{\partial{\gamma_{i_\nu}}}{\partial a_{p+M}}&
\frac{\partial{\gamma_{i_\nu}}}{\partial b_q}&\dots&\frac{\partial{\gamma_{i_\nu}}}{\partial b_{q+M}}
\end{pmatrix}
\]
As the derivative of $\gamma_{i_j}$ with respect to $b_{i_j}$ is proportional to $a_p^{i_j-q}$, the submatrix
of $D\gamma$ formed by columns corresponding to derivatives over $b_{i_1},\dots,b_{i_\nu}$ is upper triangular
with monomials $a_p^{i_1-q},\dots,a_p^{i_\nu-q}$ on the diagonal. Thus the rank of $D\gamma$ is $\nu$ as long as $a_p\neq 0$. We
conclude the proof by the implicit function theorem.
\end{proof}
\begin{remark}
Essentially the same proof can be used to show that the intersection of $L_I$ with $\{b_q=0\}$ is transverse away from $\{a_p=0\}$.
\end{remark}

It follows that $L_I$ is at least of correct codimension $\nu_I=\# I$, altough it might have singularities at the intersection with 
the hyperplane $\{a_p=0\}$.
In some cases we can explicitly describe $L_I$ and draw nice consequences.

\subsection{Explicit description of $L_I$ for $I=(q+1,\dots,q+\nu)$.}\label{sec:class}

Assume now that $I=(q+1,\dots,q+\nu)$. We can explicitly describe the set $L_I$ in that case. To avoid a flood of integer constants, let us
put
\[h_{m,n}=(q+m)-(n+p)\frac{q}{p}.\]
We will use the following simple result.
\begin{lemma}\label{lem:hisnotzero}
If $m,n\ge 0$ and $0<m+n<\frac{p+q}{\gcd(p,q)}$, then $h_{m,n}\neq 0$.
\end{lemma}
\begin{proof}

If $m=0$ then $n>0$, so $h_{m,n}\neq 0$. We can restrict ourselves to the case $m,n>0$.
We have $ph_{m,n}=pm-qn$. 
This can zero only if $m$ is an integer multiple of $\frac{q}{\gcd(p,q)}$ and $n$ is an integer
multiple of $\frac{p}{\gcd(p,q)}$.
\end{proof}
Let us also denote
\[F_{l}=\sum_{i+j=l}h_{ij}a_{p+i}b_{q+j}.\]
$F_1,\dots,F_M$ are bilinear functions on $\PP V=\PP V_x\oplus\PP V_y$.
Using this notation we can easily describe the set $L_I$.
\begin{proposition}\label{prop:LIAI}
Let $A_I=\{F_1=\dots=F_\nu\}=0$. Then $L_I\cap\{a_p\neq 0\}=A_I\cap\{a_p\neq 0\}$ as schemes, furthermore $L_I\subset A_I$. If $\nu<\frac{p+q}{\gcd(p,q)}$,
then $L_I=A_I$.
\end{proposition}
\begin{proof}
First we shall prove that $L_I$  coincides with $A_I$ away from $\{a_p=0\}$. First, a direct application of the implicit function
theorem shows that the scheme $A_I$ is smooth away from $\{a_p=0\}$. We omit the straightforward computations, which are very
similar, but simpler, to those in the proof of Lemma~\ref{lem:liissmooth}. 
Let us take an element $(a_p,\dots,a_{p+M},b_q,\dots,b_{q+M})\in L_I\cap\{a_p\neq 0\}$. The corresponding curve (see Remark~\ref{rem:elementofV})
has the following Puiseux expansion
\[y=c_qx^{q/p}+c_{q+\nu+1}x^{(q+\nu+1)/p}+\dots.\]
The coefficients $c_{q+1},\dots,c_{q+\nu}$ vanish because the parameter lies in $L_I$. By Lemma~\ref{lem:ordp1}, we infer that
if $P_1$ is the Puiseux polynomial, then
\[\ord_{t=0} P_1\ge q+\nu+p.\]
Now a straightforward computation shows that $F_j$ is the coefficient of $P_1$ at $t^{q+p+j-1}$. Therefore
\[L_I\cap\{a_p\neq 0\}\subset A_I\cap\{a_p\neq 0\}\]
as sets. The same argument shows the opposite inclusion (in the set-theoretical sense). As the two schemes are smooth, we have 
\begin{equation}\label{eq:liai}
L_I\cap\{a_p\neq 0\}=A_I\cap\{a_p\neq 0\}
\end{equation}
as schemes. Since $A_I$ is a closed scheme, $A_I\supset\ol{A_I\cap\{a_p\neq 0\}}=L_I$.

To finish the proof we need to show that if $\nu<\frac{p+q}{\gcd(p,q)}$, then the two schemes coincide on $\{a_p=0\}$. 
We shall need following fact.
\begin{lemma}\label{lem:nocomponent}
The scheme $A_I$ does not have any component lying entirely in $\{a_p=0\}$.
\end{lemma}
Given Lemma~\ref{lem:nocomponent} we finish the proof very quickly. Namely, the statement implies that
\[A_I=\ol{A_I\cap\{a_p\neq 0\}}.\]
By \eqref{eq:liai} we conclude the proof.
\end{proof}

Lemma~\ref{lem:nocomponent} shall be deduced from a more general result.

\begin{proposition}\label{prop:Nkl}
For $0\le k,l\le M$ let 
\[N_{k,l}=\{a_p=a_{p+1}=\dots=a_{p+k}=b_q=b_{q+1}=\dots=b_{q+l-1}=0\}\subset\PP V.\]
If $\nu<\frac{p+q}{\gcd(p,q)}$, then  $A_I\cap\{a_p=0\}=\bigcup_{k+l=\nu+1} N_{k,l}$ as sets. Furthermore, in the cohomology ring $H^*(\PP V;\Z)$.
\begin{equation}\label{eq:AIcapap0}
[A_I\cap\{a_p=0\}]=\sum_{k+l=\nu}\binom{\nu}{k}[N_{k,l}] .
\end{equation}
\end{proposition}
\begin{proof}
We shall proceed by induction on $\nu$. For $\nu=1$, we have $F_1=h_{10}a_{p+1}b_q+h_{01}a_pb_{q+1}$, hence $F_1\cap\{a_p=0\}$ is scheme--theoretically
$\{a_{p+1}=0\}\cup\{b_q=0\}=N_{1,0}\cup N_{0,1}$. In particular the relation~\ref{eq:AIcapap0} holds for $\nu=1$.

Suppose we have proved the statement for $\nu-1$. Let us choose $k,l\ge 0$ such that $k+l=\nu-1$. From the description of polynomials $F_1,\dots,$
we infer, that
\[F_\nu|_{N_{k,l}}=h_{k+1,l+1}a_{k+1}b_{l+1}\neq 0.\]
Hence,
\[N_{k,l}\cap\{F_\nu=0\}=N_{k+1,l}\cup N_{k,l+1}\]
In particular
\[[N_{k,l}\cap\{F_\nu=0\}]=[N_{k+1,l}]+[N_{k,l+1}]\]
in the cohomology ring. The induction step follows from the well-known formula $\binom{a}{b}=\binom{a-1}{b-1}+\binom{a-1}{b}$.
\end{proof}

\begin{proof}[Proof of Lemma~\ref{lem:nocomponent}]
As $A_I$ is given by $\nu$ equations, each component of $A_I$ must have codimension at most $\nu$ (see e.g. \cite[Theorem~I.7.2]{Har}). 
On the other hand, all components
of $A_I\cap\{a_p=0\}$ have codimension $\nu+1$.
\end{proof}

\subsection{Proof of Theorem~\ref{thm:converse}}\label{sec:converse}
Let us recall the formulation of Theorem~\ref{thm:converse}.

\begin{theorem}
Let $p,r_0,r_1,p'$ be positive integers such that $r_1>r_0$, $p'>p$, $r_0<r_1+(p-p')$ and
\[r_1-r_0<\frac{p+r_0}{\gcd(p,r_0)}.\] 
Then, given any Puiseux expansion
\begin{equation}\label{eq:referencePuis2}
y=d_{r_1+p-p'}x^{(r_1+p-p')/p'}+d_{r_1+p-p'+1}x^{(r_1+p-p'+1)/p'}+\dots,
\end{equation}
and for an arbitrary $N\ge 0$, there exists a deformation such that for $s\neq 0$
for $s\neq 0$, the corresponding singularity has the expansion starting with
\[y_s=c_{r_0}x_s^{r_0/p}+c_{r_1}x_s^{r_1/p}+\dots,\]
where $c_{r_0}c_{r_1}\neq 0$, and
for $s=0$, the first $N$ terms of the Puiseux expansion agree with those of \eqref{eq:referencePuis2}.
\end{theorem}
\begin{proof}
Let $q=r_0$.
Let us choose $M$ (used in the definition of $\PV$) equal to $N+r_0+p'-p$. Let us choose a point $(x',y')\in\PP V$ such that
$a_p=\dots=a_{p'-1}=0$, $a_{p'}=1$, $a_{p'+1}=a_{p'+2}=\dots=0$ and $b_q=\dots=b_{r_1+p-p'-1}=0$ and for all $j\ge r_1+p-p'$ we have
$b_j=d_j$, where $d_j$ comes from \eqref{eq:referencePuis2}. It is clear, that to the point $(x',y')$ corresponds a curve
which has Puiseux expansion agreeing with \eqref{eq:referencePuis2} up to the term $x^{(q+M)/p'}$.

The point $(x',y')$ lies in $N_{k,l}$ for $k=p'-p$ and $l=r_1+p-p'-r_0$ (see Proposition~\ref{prop:Nkl}). 
Let $\nu=k+l$ and $I=(q+1,\dots,q+\nu)$ and $I'=(q+1,\dots,q+\nu+1)$.  We have $\nu<\frac{p+r_0}{\gcd(p,r_0)}$ by assumptions of the theorem. By 
Proposition~\ref{prop:Nkl} we infer that $N_{k,l}\subset A_I\cap\{a_p=0\}$. Let us now choose a generic hypersurface $P$ of codimension $\nu+1$, which
passes through $(x',y')$ and intersects $A_{I'}$ in finitely many points. Such $P$ exists, because $A_{I'}$ has codimension $\nu+1$. Indeed, if 
$\nu+1<\frac{p+r_0}{\gcd(p,r_0)}$, then $\codim A_{I'}=\nu+1$ by Proposition~\ref{prop:Nkl}. If $\nu+1=\frac{p+r_0}{\gcd(p,r_0)}$, then 
$A_{I'}\cap\{a_p\neq 0\}$ has codimension $\nu+1$, and $A_{I'}\cap\{a_p=0\}\subset A_{I}\cap\{a_p=0\}$, which has also codimension $\nu+1$ in $\PV$. In
other words, $A_{I'}$ might have components in $\{a_p=0\}$, but no component of improper codimension.

The intersection $P\cap A_{I}$ is a complex curve, smooth away from $\{a_p=0\}$.
Let us choose a neighborhood of $U$ of $(x',y')$ in $P\cap A_{I}$. $U$ can be chosen small enough so that $U\cap A_{I'}=(x',y')$ (because
the intersection of $P$ with $A_{I'}$ is finite). Shrinking $U$ if necessary we can furthermore assume that $U\cap\{a_p=0\}=(x',y')$. Furthermore,
$U$ is smooth away from $(x',y')$. Now the normalization of $U$, denoted by $\widehat{U}$, is a smooth curve. $z\in\widehat{U}$ be a preimage
of $(x',y')$ under the normalization map. Let $D$ be a neighborhood of $z$. Then $D$ provides a required deformation.

More precisely, let $\pi\colon\widehat{U}\to U\subset\PP V$ be the normalization map. Let us choose a local variable $s$ on $D$ such that $s=0$ corresponds to $z$.
The map $\pi$ restricted to $D$ can be written as
\[s\to(a_p(s),a_{p+1}(s),\dots,a_{p+M}(s),b_q(s),\dots,b_{q+M}(s))\in\PV.\]
We can lift $\pi$ to a map $\tilde{\pi}\colon\widehat{U}\to V$.
The lift $\tilde{\pi}$ is not unique, but we can
choose the one that $a_{p'}(0)=1$, $b_j(0)=d_j$ for all $j\ge r_1+p-p'$.
All the required properties from the statement of the theorem are satisfied by construction. 
\end{proof}

\begin{example}\label{ex:1016again}
Let $p=9$, $p'=10$, $r_0=9$ and $r_1=17$. All the assumptions of Theorem~\ref{thm:converse} are satisfied. As $r_1+p-p'=16$, there
exists a deformation with $y_s=c_{9}(s)x_s^{9/9}+c_{17}(s)x_s^{17/9}+\dots$ and $y_0=x_0^{16/10}+x_0^{17/10}+\dots$.
\end{example}

\subsection{Cohomology classes of $L_I$}\label{sec:rhisone}
We now revert to a general $I$. Let $\nu=\nu_I=\# I$.
Let us now denote by $H_x$ (respectively $H_y$) the hyperplane class in $H^2(\PP V_x;\Z)$ (respectively $H^2(\PP V_y;\Z)$). The cohomology ring of $\PP V$ is
\[H^*(\PP V;\Z)=\Z[H_x]\oplus\Z[H_y]/(H_x^M-1,H_y^M-1).\]
The cohomology class of the scheme $L_I$ in $H^{2\nu}(\PP V;\Z)$ can be written as
\[[L_I]=\sum_{k=0}^{\nu_I} l_k\cdot H_x^kH_y^{\nu-k},\]
for some numbers $l_0,\dots,l_{\nu}$. The integers $l_0,\dots,l_{\nu}$ are in general hard to compute. In special cases we can do that.
\begin{example}
Assume that $I=(q+1,\dots,q+\nu)$ and $\nu<\frac{p+q}{\gcd(p,q)}$. Then $l_k=\binom{\nu}{k}$. This follows directly from Proposition~\ref{prop:LIAI},
since the class of $N_{k,l}$ is $H_x^{k+1}H_y^l$.
\end{example}

In the general case we have the following result.

\smallskip
\begin{proposition}~\label{prop:rhisone}
For any $I$, we have $l_{\nu}=1$.
\end{proposition}
\begin{proof}
Let us intersect $L_I$ with a generic plane $P$ in the class $H_x^MH_y^{M-\nu}$. This intersection is
a finite set of points. As $L_I$ has no component lying entirely in $\{a_p=0\}$ (and also no component lying entirely in $\{b_q=0\}$), 
by picking a sufficiently generic $P$ we may assume that no
intersection point lies on $\{a_pb_q=0\}$. We want to show that there is only one point in the intersection.

Let us pick affine coordinates on $\PP V\setminus\{a_pb_q=0\}$ (which still we denote by
$a_{p+1},\dots,a_{p+M}$,$b_{q+1},\dots,b_{q+M}$). Then, the plane $P$ is given by $M$ equations
of the form
\begin{align}
\theta_1a_{p+1}+\dots+\theta_Ma_{p+M}&=\theta_0\label{eq:linearina}\\
\intertext{and $M-\nu$ equations of the form}
\theta'_1b_{q+1}+\dots+\theta'_db_{q+M}&=\theta'_0,\label{eq:linearinb}
\end{align}
where $\theta_0,\dots,\theta_M,\theta'_0,\dots,\theta'_M$ are generic (of course we take a different set of $\theta$'s for different equations).
Let $I=(i_1,\dots,i_\nu)$. Then $L_I$ in the affine part can be presented by
\begin{equation}\label{eq:Pi's}
\begin{split}
\gamma_{i_1}(1,a_{p+1},\dots,a_{i_1},1,b_{q+1},\dots,b_{i_1})&=0\\
\dots\\
\gamma_{i_\nu}(1,a_{p+1},\dots,a_{i_\nu},1,b_{q+1},\dots,b_{i_\nu})&=0.
\end{split}
\end{equation}

We want to show that \eqref{eq:linearina}, \eqref{eq:linearinb} and \eqref{eq:Pi's} have
a unique solution. The equations from \eqref{eq:linearina} uniquely determine $a_{p+1},\dots,a_{p+M}$.
Then \eqref{eq:linearinb} and \eqref{eq:Pi's} become a system of linear equations on $b$'s. We need to show that
the this system is non--degenerate. By genericity of $P$, it is enough to show that the matrix of coefficients
of \eqref{eq:Pi's} has rank $\nu$ (the $(j,k)$ entry of this matrix is the coefficient at $b_{q+k}$ of the $j$-th equation, i.e. at $\gamma_{i_j}$).
The argument resembles proof of Lemma~\ref{lem:liissmooth}.

The submatrix formed by $(j,i_k)$ entries when $1\le j,k\le\nu$ is lower triangular with $1$'s at the diagonal.
Indeed, if $k>j$, then the coefficient at $b_{q+k}$ in  $\gamma_{i_j}$ is zero. The coefficient at $b_{q+k}$ in $\gamma_{i_k}$ is
equal to $1$, this follows from the proof of Proposition~\ref{prop:homog}.
\end{proof}

\subsection{Computation of $l_0$.}\label{sec:whatisl0}
In this section we shall prove, that under some additional assumptions on $I$, we have $l_0=1$ as well.
The trick is to consider the reverse Puiseux expansion
\begin{equation}\label{eq:revpuis}
x=g_{p}y^{p/q}+g_{p+1}y^{(p+1)/q}+\dots.
\end{equation}
If we are given the expansions \eqref{eq:puis} and \eqref{eq:revpuis} we can pass from one to another, i.e. express $g_p,\dots$ as functions
of $c_q,c_{q+1},\dots$. We have the following simple result.

\begin{lemma}\label{lem:gishomog}
For any $i>0$ the function $(g_p,g_{p+1},\dots,)\to c_{q+i}\cdot g_p^{(q+i)/p+i}$ is a homogeneous polynomial of degree $i$, which is also weighted homogeneous
of degree $i$ if the weight of $g_{p+j}$ is defined as $j$. Conversely, for any $j>0$, the function $(c_q,c_{q+1},\dots,)\to g_{p+j}\cdot c_q^{(p+j)/q+j}$
are homogeneous polynomials of weight $j$ and also weighted homogeneous of degree $j$ if the weight of $c_{q+i}$ is defined as $i$.
\end{lemma}
\begin{proof}
We shall prove only the first part, the other is completely analogous.
Let us denote $\tau=y^{1/q}$. By \eqref{eq:revpuis}, the singularity can be parametrized as
\begin{align*}
x&=g_p\tau^p+g_{p+1}\tau^{p+1}+\dots\\
y&=\tau^q.
\end{align*}
The statement follows now immediately from Proposition~\ref{prop:homog}.
\end{proof}

Let $U_1,U_2,\dots$ be polynomials defined by the property that
\begin{equation}\label{eq:uj}
U_j(c_q,c_{q+1},\dots,c_{q+j})=g_{p+j}\cdot c_q^{(q+i)/p+i}.
\end{equation}

Similarly as in \eqref{eq:cqi} let us define functions $\rho_{p},\rho_{p+1},\dots$ by the formula

\begin{equation}\label{eq:rhopi}
g_{p+j}=b_q^{-(p+j)/p-j}\rho_{p+j}.
\end{equation}

\begin{lemma}
For any $j\ge 0$ we have
\[\rho_{p+j}=U_j(\gamma_{q},\gamma_{q+1},\dots,\gamma_{q+j})a_p^{1-j}.\]
\end{lemma}
\begin{proof}
Substituting \eqref{eq:cqi} into $U_j$ and using the homogeneity of $U_j$ we obtain
\[a_p^{-qj/p-j/p-j}U_j(\gamma_q,\dots,\gamma_{q+j})=g_{p+j}c_q^{-(q+j)/p-j}.\]
By \eqref{eq:rhopi} we infer that
\[U_j(\gamma_q,\dots,\gamma_{q+j})a_p^{-qj/p-j/p-j}b_q^{(p+j)/p+j}c_q^{-(p+j)/q-j}=g_{p+j}.\]
But $c_q=b_qa_p^{-q/p}$. The lemma follows.
\end{proof}

\begin{example}\label{ex:gone}
Changing roles of $a$'s and $b$'s in \eqref{eq:stupidgamma} yields $\rho_{p+1}=a_{p+1}b_q-\frac{p}{q}a_pb_{q+1}$. Thus $\rho_{p+1}=-\frac{p}{q}\gamma_{q+1}$.
\end{example}

For the set of indices $I'$ we can define a subscheme of $\PP V$:
\[R_{I'}=\ol{\{(a_p,\dots,b_q,\dots)\colon\forall_{j\in I'}\rho_j=0\}\setminus\{b_q=0\}}.\]
by a complete analogy to Definition~\ref{def:LI}. By Lemma~\ref{lem:liissmooth}, $R_{I'}$ is smooth away from $\{b_q=0\}$. 

For any $I$ let us denote
\[\Gamma_I=\{i\in\mathbb{N}\colon i+q\not\in I\}.\]
If $\kappa=(p;q_1,\dots,q_n)$, $q=p$ and $I=I_\kappa$, by Lemma~\ref{lem:semigroup} we know that $\Gamma_I$ is a semigroup. 
If $q=q_1$ then the argument from proof of Lemma~\ref{lem:semigroup}
can be used to show that $\Gamma_{I^q_\kappa}$ is a semigroup as well. These two cases are the most interesting ones.

\begin{proposition}\label{prop:reverse}
Assume that the set $\Gamma_I$ is a semigroup. Let $I'=(j\colon i+q-p\in I)$. Then $L_I$ and $R_{I'}$ coincide.
\end{proposition}
\begin{proof}
By induction on $\#I$ we shall show that $L_I$ and $R_{I'}$ coincide away of $\{a_p=0\}\cup\{b_q=0\}$. 
Let $I=(i_1,\dots,i_\nu)$. We shall denote by $I_k$ the subset $(i_1,\dots,i_k)$ consisting on first $k$ elements of $I$, 
and $I'_k=(j\colon i+q-p\in I)$. Obviously, for any $k$
the set $\{i\in\mathbb{N}\colon i+q\not\in I_k\}$ is still a semigroup.

For $k=1$, the semigroup condition means that $I_1=(q+1)$. By Example~\ref{ex:gone} we have the coincidence of $L_{I_1}$ and $R_{I'_1}$. Assume that we proved that
$L_{I_{k-1}}$ coincides with $R_{I_{k-1}}$ away from $\{a_p=0\}\cup\{b_q=0\}$. Observe that if $a_p,b_q\neq 0$, we have $\gamma_p\neq 0$ and $\rho_q\neq 0$.

Let us take $i_k$. For any presentation
\[i_k-q=\sum_{s=1}^h (m_s-q),\]
i.e. for any choice of finitely many integers $m_1,\dots,m_h\ge q$ 
satisfying the above condition, there must exists at least index $s_0$ such that $m_{s_0}\in I_k$, otherwise
the semigroup assumption is violated. This means that any monomial of the form $\gamma_{q+1}^{k_1}\cdot \dots\cdot\gamma_{q+M}^{k_M}$ such that 
$\sum ik_i=i_k-q$
must vanish on $L_{I_k}$, because for some $i$ such that $i\in I_{k}$ we have $k_i>0$, and then $\gamma_{q+i}^{k_i}$ vanishes on $L_{I_k}$.
 
By the homogeneity of $U_{i_k}$ it follows that $U_{i_k}$ vanishes on $L_{I_k}$, 
hence $\rho_{i_k+p-q}$ vanishes on $L_{I_k}$. In particular $L_{I_k}\subset R_{I_k}$
away from $\{a_p=0\}\cup\{b_q=0\}$. The opposite inclusion is proven in the same way. The induction step is finished.

We clearly have $L_I=\ol{L_I\setminus(\{a_p=0\}\cup\{b_q=0\})}$. Indeed, $L_I=\ol{L_I\setminus\{a_p=0\}}$ and $L_I\setminus\{a_p=0\}$ is smooth and not contained
in $b_q=0$, so $\ol{L_I\setminus(\{a_p=0\}\cup\{b_q=0\})}=\ol{L_I\setminus\{a_p=0\}}=L_I$. Similar identity holds for $R_{I'}$
We conclude that $L_I=R_{I'}$.
\end{proof}
The result shows, in particular, that the only singularities of $L_I$ can occur on $\{a_p=b_q=0\}$, because $R_{I'}$ is smooth away from $\{b_q=0\}$.
We are ready to state the result about $l_0$. We call it a theorem, because of the importance in the next section.

\begin{theorem}\label{thm:l0isone}
Assume that $\Gamma_I$ is a semigroup. Then $l_0=1$.
\end{theorem}

\begin{proof}
By Proposition~\ref{prop:reverse} we have $L_I=R_{I'}$. The argument of Proposition~\ref{prop:rhisone} applied to $R_{I'}$ 
shows that $l_0=1$.
\end{proof}

The assumption on $\Gamma_I$ being a semigroup is important, because of the following annoyingly simple counterexample.

\begin{example}
Assume that $I=(q+j)$. Then $l_0=j$. Indeed, the scheme $\{\gamma_{q+j}=0\}$ has no component lying entirely in $\{a_p=0\}$, because $a_p$ does
not divide $\gamma_{q+j}$. Hence $L_I=\{\gamma_{q+j}=0\}$. By Proposition~\ref{prop:homog}, $[L_I]=jH_x+H_y$.
\end{example}

Theorem~\ref{thm:l0isone} is somewhat unexpected if one looks at the degrees of the functions $\gamma_{i_1},\dots,\gamma_{i_\nu}$ 
with respect to the
variables $a_p,\dots,a_{p+M}$. A naive degree counting argument suggests that $l_0=(i_1-q)\cdot(i_2-q)\cdot\ldots\cdot(i_\nu-q)$. 
This is very far from true and suggests a complicated behavior
of the polynomials $\gamma_{q+1},\gamma_{q+2},\dots$.

\subsection{Non--genericity of the functions $\gamma_i$}\label{sec:bernkou}

For $i>0$, let $N_i$ be the Newton polytope related to the polynomial $\gamma_{q+i}$. The assumption, that for each $i>0$,
the function $\gamma_{q+i}$ is generic in the space of polynomials having $N_i$ as their Newton polytope, will lead
to a wrong prediction of the number $l_0$ for some sets $I$. To simplify the argument we
deal mainly with the case $I=(q+1,q+3,q+5)$, however the methods we use can be generalized.

\begin{proposition}\label{prop:135}
Let $f_1,f_3,f_5$ be generic polynomials with Newton polytope $N_1,N_3$ and $N_5$. Then, for generic plane $P\subset\PP V$ representing
the class $H_x^{M-3}H_y^M$, there exist two points of intersection of $P$ with the set $\{f_1=f_3=f_5=0\}$, which lie away from the 
hypersurfaces $\{a_p=0\}$, $\{b_q=0\}$.
\end{proposition}
\begin{proof}
To study intersections on $\PP V\setminus\{a_pb_q=0\}$, we choose charts on $\PP V_x\setminus\{a_p=0\}$ and $\PP V_y\setminus\{b_q=0\}$.
We shall still call the coordinates $a_{p+1},\dots,a_{p+M}$ and $b_{q+1},\dots,b_{q+M}$. For simplicity, let us assume that $M=5$.
The plane $P$ is given by the equations
\begin{equation}\label{eq:Pfirst}
\begin{split}
\alpha_{11} a_{p+1}+\alpha_{12} a_{p+2}+\alpha_{13} a_{p+3}+\alpha_{14} a_{p+4}+\alpha_{15}a_{p+5}&=\alpha_{10}.\\
\alpha_{21} a_{p+1}+\alpha_{22} a_{p+2}+\alpha_{23} a_{p+3}+\alpha_{24} a_{p+4}+\alpha_{25}a_{p+5}&=\alpha_{20}.\\
\beta_{j1} b_{q+1}+\beta_{j2} b_{q+2}+\beta_{j3} b_{q+3}+\beta_{j4} b_{q+4}+\beta_{j5}b_{q+5}&=\beta_{j0}.
\end{split}
\end{equation}
The last equation is repeated 5 times for $j=1,\dots,5$. The complex coefficients 
$\alpha_{ij}$ and $\beta_{ij}$ are generic. Generic $f_1,f_3,f_5$ have form (it should be understood that before each monomial
in the following equation stays a generic coefficient):
\begin{equation}\label{eq:Ffirst}
\begin{split}
0=&b_{q+1}+a_{p+1}\\
0=&b_{q+3}+b_{q+2}a_{p+1}+b_{q+1}a_{p+1}^2+a_{p+1}^3+b_{q+1}a_{p+2}+a_{p+1}a_{p+2}+a_{p+3}\\
0=&b_{q+5}+b_{q+4}a_{p+1}+b_{q+3}a_{p+1}^2+b_{q+2}a_{p+1}^3+b_{q+1}a_{p+1}^4+a_{p+1}^5+\\
&+b_{q+3}a_{p+2}
+b_{q+1}a_{p+2}^2+b_{q+2}a_{p+3}+b_{q+1}a_{p+4}+a_{p+5}+a_{p+2}a_{p+3}.
\end{split}
\end{equation}
The ten equations will be solved in the following way.
\begin{itemize}
\item[$\bullet$] From last five equations of \eqref{eq:Pfirst} (i.e. the the third equation repeated 5 times), we
determine uniquely $b_{q+1},\dots,b_{q+5}$.
\item[$\bullet$] The first equation of \eqref{eq:Ffirst} give us $a_{p+1}$.
\item[$\bullet$] The second equation of \eqref{eq:Ffirst} lets us express $a_{p+3}$ as a linear function of $a_{p+2}$.
\item[$\bullet$] From the first two equations of \eqref{eq:Pfirst} we express $a_{p+4}$ and $a_{p+5}$ as linear, non--homogeneous, functions of $a_{p+2}$.
\item[$\bullet$] After all substitutions, the last equation of \eqref{eq:Ffirst} takes form
\[
\alpha_0 a_{p+2}^2+\alpha_1a_{p+2}+\alpha_2=0,
\]
where by choosing suitable coefficients of $f_1,f_2,f_3$ and coefficients entering in $P$ we have a full control
over the coefficients $\alpha_0,\alpha_1$ and $\alpha_2$. In particular, we see, that for generic coefficients, the equation has two
solutions, which correspond to the two distinct intersection points of $\{f_1=f_3=f_5=0\}$ with $P$.
\end{itemize}
\end{proof}

In general, for other sets $I$, the coefficient $l_0$ --- assuming genericity of $\gamma_{q+i}$'s --- can
be computed using techniques of Kuschnirenko and Bernstein (see \cite[Section 5.5]{Fu} or \cite{Ku,Be}). Their method is as follows. Consider
$N_A$ and $N_B$ the Newton polygon related to the linear polynomial $\theta_1a_{p+1}+\dots+\theta_Ma_{p+M}+\theta_0$ (respectively,
$\theta_1b_{q+1}+\dots+\theta_Mb_{q+M}+\theta_0$), for $\theta_0,\dots,\theta_M$ non--zero. For $N_i^0$ be the Newton polygon of
$\gamma_{q+i}$ restricted to $\{a_p=b_q=1\}$. The number of non--zero solutions 
to the system 
\begin{align*}
\gamma_{i_1}=\dots=\gamma_{i_\nu}&=0\\
\alpha_{i1} a_{p+1}+\dots+\alpha_{jM}a_{p+M}&=\alpha_{i0}\textrm{ ($M-\nu$ equations for $i=1,\dots,M-\nu$)}\\
\beta_{j1} b_{q+1}+\dots+\beta_{jM}b_{q+M}&=\beta_{j0}\textrm{ ($M$ equations for $j=1,\dots,M$)},
\end{align*}
where $\alpha_{ik}$, $\beta_{jk}$ ($1\le i\le M-\nu$, $1\le j\le M$, $0\le k\le M$) are generic,
is bounded from above by
\begin{equation}\label{eq:mvol}
(2M)!\cdot\MVol(N_{i_1}^0,\dots,N_{i_\nu}^0,\underbrace{N_A,\dots,N_A}_{\text{$M-\nu$ times}}
\underbrace{N_B,\dots,N_B}_{\text{$M$ times}}).
\end{equation}
Here $\MVol$ denotes the so-called mixed volume of a system of polytopes. Furthermore, under explicit genericity condition
(see \cite[Theorem B]{Be}, \cite[Paragraph 1.19.III']{Ku} 
or \cite[Section 5.5]{Fu}), the number of solution is actually equal to \eqref{eq:mvol}. 
In particular, whenever the value \eqref{eq:mvol} is greater than 1, 
we know that these genericity conditions are violated.
Proposition~\ref{prop:135} implies that these genericity
conditions are not satisfied if $I=(q+1,q+3,q+5)$. We do knot know a formula for the mixed volume for general $I$. For specific values, it might be
often computed using a computer program.

\section{WDVV equations}\label{S:WDVV}
The failure of genericity of Puiseux coefficients described in Sections~\ref{sec:rhisone} and~\ref{sec:bernkou} was detected
by looking at the Puiseux expansions of $y$ in powers of $x$ and the Puiseux expansion of $x$ in powers of $y$. This way of thinking
can be continued. We shall show that functions $\gamma_{q+1},\gamma_{q+2},\dots$ (regarded as functions on the space of parameters)
have various symmetries. These symmetries can be encoded in WDVV equations.

\subsection{Review of WDVV equations}\label{sec:reviewWDVV}
Notes by Dubrovin \cite{Du} are an excellent introduction for WDVV equations. Here we review only necessary basics.
Consider a $C^3$ function $F(t_1,\dots,t_n)$, where $t_1,\dots,t_n$ are variables (real or complex, but in
the latter case we require $F$ to be holomorphic). For $\alpha,\beta,\gamma=1,\dots,n$, let us denote
\[c_{\alpha\beta\gamma}=\frac{\p^3 F}{\p t_\alpha\p t_\beta\p t_\gamma}.\]
We write
\[\eta_{\beta\gamma}=c_{1\beta\gamma}.\]
We shall treat $\eta_{\beta\gamma}$ as a metric on $\R^n$. More precisely, let $H=\{\eta_{\beta\gamma}\}_{\beta,\gamma=1}^n$. 
Assume that $H$ is non--degenerate and write $\eta^{\beta\gamma}=(H^{-1})_{\beta,\gamma}$. We define
\[c_{\alpha\beta\gamma}=\sum_{\delta=1}^n\eta^{\alpha\delta}c_{\delta\beta\gamma}.\]
Let us choose a basis $e_1,\dots,e_n$ of $\R^n$. We define a scalar product in that basis by
\[\langle e_\alpha,e_\beta\rangle=\eta_{\alpha\beta}.\]
The coefficients $c_{\alpha\beta\gamma}$ are used to define a multiplication $\R^n\times\R^n\to\R^n$ by the formula
\[e_\beta\cdot e_\gamma=\sum_{\alpha,\delta=1}^n\eta^{\alpha\delta}c_{\delta\beta\gamma}e_\alpha.\]
As the coefficients $c_{\alpha\beta\gamma}$ depend in general on $t_1,\dots,t_n$, strictly speaking we have a family of multiplications. Each
such multiplication is commutative (by the symmetry of third derivatives) but not always associative. The associativity holds if
and only if for any $\alpha,\beta,\gamma,\delta$ the following equation is satisfied.
\begin{equation}\label{eq:WDVV-standard}
\frac{\p^3 F}{\p t_\alpha\p t_\beta\p t_\lambda}\eta^{\lambda\mu}\frac{\p^3 F}{\p t_\gamma\p t_\delta\p t_\mu}=
\frac{\p^3 F}{\p t_\gamma\p t_\beta\p t_\lambda}\eta^{\lambda\mu}\frac{\p^3 F}{\p t_\alpha\p t_\delta\p t_\mu},
\end{equation}
where we sum over repeated indices $\lambda$ and $\mu$.

\begin{definition}
The equations \eqref{eq:WDVV-standard} is called the \emph{WDVV equations}. A function $F$ such that the matrix $H$ defined above is non--degenerate
and which satisfies \eqref{eq:WDVV-standard} is called the \emph{WDVV potential}.
\end{definition}

Given the complexity of the WDVV equations it is hard to construct non--trivial solutions (apart from the obvious ones, 
where the third derivative is constant). Highly non--trivial solutions appear in many contexts, from
topological quantum field theories, through Gromov--Witten potentials up to KdV hierarchies. 
Even a short review of these application is beyond the range of this article. We point out that K.~Saito \cite{SK1,SK2} and M.~Saito \cite{SM}
found a WDVV potential (more precisely, the structure of a Frobenius manifold, which is tightly related to WDVV equations, but
we do not discuss it here)
related with deformation of the hypersurface singularities (see \cite[Chapter 11nn]{He}). The underlying
vector space is the tangent space to the base of semiuniversal unfolding and its dimension is equal to the Milnor number of a corresponding singularity.
We do not know, how this structure can be related to our construction, that we describe below.

\subsection{Introducing new notation}
In Section~\ref{S:WDVV} we shall assume that the singularity has the following Puiseux expansion
\begin{equation}\label{eq:yofx}
y=x^{q/p}+c_1x^{(q+1)/p}+c_2x^{(q+2)/p}+\dots
\end{equation}
In particular, the coefficient at $x^{q/p}$ is always assumed to be $1$. Furthermore we shift the indices, otherwise
the formulas would be very difficult to read. The reverse Puiseux expansion will be defined as
\begin{equation}\label{eq:xofy}
x=y^{p/q}+g_1y^{(p+1)/q}+g_2y^{(p+2)/q}+\dots
\end{equation}
There is a relation between the coefficients $c_1,\dots$ and $g_1,\dots$. Namely, when we replace $x$ in \eqref{eq:yofx} by its the expansion \eqref{eq:xofy},
we should get $y=y$. More precisely we have the following result.
\begin{lemma}
There exist polynomials with rational coefficients $W_1,W_2,\dots,$ such that $W_i$ depends only on $c_1,\dots,c_i$ 
such that for each $i=1,\dots$
we have
\[g_i=W_i(c_1,\dots,c_{i}).\]
Furthermore, $W_i$ is weighted homogeneous of weight $i$, if the weights of $c_1,c_2,\dots$ are, 
respectively $1,2,\dots$.
\end{lemma}
\begin{proof}
This follows immediately from Lemma~\ref{lem:gishomog}, because 
\[W_i(x_1,\dots,x_i)=U_i(1,x_1,\dots,x_i),\]
where $U_i$ is as defined in \eqref{eq:uj} (remember that in this section we shift indices, so $c_{q+i}$ in previous sections corresponds
to $c_i$ now).
\end{proof}
\begin{remark}
Unlike $U_i$, the polynomial $W_i$ is no longer homogeneous itself (unlike  Lemma~\ref{lem:gishomog}), 
because we normalized $c_0=1$.
\end{remark}

We are going now to state the following result, which will show the symmetry of functions $g_1,g_2,\dots$.

\begin{theorem}\label{thm:derg}
The derivative 
$\frac{\partial^r g_n}{\partial c_1^{r_1}\dots\partial c_n^{r_n}}$
depends only on $n$, $r$ and $\sum ir_i$. In other words, if we are given two sets of numbers $r_1,\dots,r_n$ and $r'_1,\dots,r_n'$
with $\sum r_i=\sum r_i'$ and $\sum ir_i=\sum ir_i'$, then
\[\frac{\partial^r g_n}{\partial c_1^{r_1}\dots\partial c_n^{r_n}}=
\frac{\partial^r g_n}{\partial c_1^{r'_1}\dots\partial c_n^{r'_n}}.\]
\end{theorem}
The proof of Theorem~\ref{thm:derg} is contained in Section~\ref{sec:derg}. Now we shall introduce some additional notation. Let
\begin{align*}
u&=x^{1/p}\\
w&=y^{1/q}\\
z&=(1+g_1 w+\dots)^{1/p}\\
G&=g_1 w+g_2 w^2+\dots\\
G_i&=\frac{\partial g_1}{\partial c_i}w+\frac{\partial g_2}{\partial c_i}w^2+\dots\\
G_{ij}&=\frac{\partial^2 g_1}{\partial c_i\partial c_j}w+\frac{\partial^2 g_2}{\partial c_i\partial c_j}w^2+\dots\\
G_{ijk}&=\frac{\partial^3 g_1}{\partial c_i\partial c_j\partial c_k}w+\frac{\partial^3 g_2}{\partial c_i\partial c_j\partial c_k}w^2+\dots.
\end{align*}
Here and afterwards, $i,j,k$ are indices from $1$ to infinity.
The series $G$ can be regarded as a generating function for $g_1,g_2,\dots$. We shall always treat the variables $g_1,g_2,\dots$ as functions
of $c_1,c_2,\dots$ and $u$.

\subsection{Proof of Theorem~\ref{thm:derg}}\label{sec:derg}
We begin with computing $G_{ij}$ and showing that it depends only on $i+j$. Our first step is computing $G_i$. Let us first write
\begin{equation}\label{eq:xyyx}
\begin{split}
w^q=&u^q+c_1u^{q+1}+c_2u^{q+2}+\dots\\
u^p=&w^p+g_1w^{p+1}+g_2w^{p+2}+\dots
\end{split}
\end{equation}
Differentiating the second equation of \eqref{eq:xyyx} over $c_i$ we obtain
\[
0=\frac{\partial w^p}{\partial c_i}+g_1\frac{\partial w^{p+1}}{\partial c_i}+\frac{\partial g_1}{\partial c_i}w^{p+1}+
g_2\frac{\partial w^{p+2}}{\partial c_i}+\frac{\partial g_2}{\partial c_i}w^{p+2}+\dots
\]
Thus
\[-w^p\left(\frac{\partial g_1}{\partial c_i}w+\frac{\partial g_2}{\partial c_i}w^2+\dots\right)=%
\left(pw^{p-1}+g_1(p+1)w^p+\dots\right)\cdot\frac{\partial w}{\partial c_i}.\]
But
\[
pw^{p-1}+g_1(p+1)w+\dots=\frac{\partial}{\partial w}(w^pz^p)=pw^{p-1}z^p+pw^pz^{p-1}\frac{\partial z}{\partial w}.
\]
Hence
\begin{equation}\label{eq:G1}
-w^pG_i=pw^{p-1}z^{p-1}\left(z+w\frac{\partial z}{\partial w}\right)\frac{\partial w}{\partial c_i}.
\end{equation}
We use now that $w=(u^q+c_1 u^{q+1}+\dots)^{1/q}=u(1+c_1u+\dots)^{1/q}$.
Hence
\begin{equation}\label{eq:zz2}
\frac{\partial w}{\partial c_i}=\frac{1}{q}u^{i+1}(1+c_1u+\dots)^{1/q-1}=\frac{1}{q}w^{1-q}u^{q+i}.
\end{equation}
Recalling that $u=wz$ we obtain
\[
\frac{\partial w}{\partial c_i}=\frac{1}{q}w^{i+1}z^{q+i}
\]
Finally we obtain the desired formula for $G_i$.
\begin{equation}\label{eq:Gi}
G_i=-\frac{p}{q}w^{i}z^{p+q+i-1}\left(z+w\frac{\partial z}{\partial w}\right).
\end{equation}

In order to compute $G_{ij}$ we differentiate \eqref{eq:Gi} over $c_j$. We get
\begin{multline*}
\frac{\partial G_i}{\partial c_k}=-\frac{p}{q}w^{i-1}z^{p+q+i-2}\left[%
\left((i+1)wz\frac{\partial z}{\partial w}+iz^2\right)\frac{\partial w}{\partial c_k}\right.\\
+\left.\left((p+q+i)wz+(p+q+i-1)w^2\frac{\partial z}{\partial w}\right)\frac{\partial z}{\partial c_k}+
w^2z\frac{\partial^2z}{\partial c_k\partial w}\right].
\end{multline*}
To compute $\frac{\partial z}{\partial c_k}$ observe that similarly as in \eqref{eq:zz2} we have
\[
\frac{\partial z}{\partial g_l}=\frac{1}{p}w^lz^{1-p}
\]
Hence, by the chain rule
\[
\frac{\partial z}{\partial c_k}=
\sum_{l=1}^\infty\frac{\partial z}{\partial g_l}\frac{\partial g_l}{\partial c_k}=\frac{1}{p}z^{1-p}\sum_{l=1}^\infty w^{l}
\frac{\partial g_l}{\partial c_k}=-\frac{1}{q}w^kz^{q+k}\left(w+z\frac{\partial z}{\partial w}\right).
\]
Differentiating that with respect to $w$ we obtain
\begin{align*}
\frac{\partial^2 z}{\partial c_k\partial w}&=
-\frac{1}{q}w^{k-1}z^{q+k-1}\cdot\\\cdot&\left[kz^2+(q+2+2k)wz\frac{\partial z}{\partial w}+(q+k)w^2\left(\frac{%
\partial z}{\partial w}\right)^2+w^2z\frac{\partial^2z}{\partial w^2}\right].
\end{align*}
Then we get the following result
\begin{equation}\label{eq:zz12}
\begin{split}
G_{ik}&=\frac{p}{q}w^{i+k}z^{p+2q-2+i+k}\cdot\\&\cdot\left[(p+q+i+k)z^2+(2p+3q+2i+2k+1)zw\frac{\partial z}{\partial w}+\right.\\
&+\left.%
(p+2q+i+k-1)w^2\left(\frac{\partial z}{\partial w}\right)^2+w^2z\frac{\partial^2 z}{\partial w^2}\right]
\end{split}
\end{equation}
We can see that $G_{ik}$ depends on $i$ and $k$ only through the sum $i+k$. It is straightforward to show that $G_{ijk}$ depends only on $i+j+k$.
We have
\begin{multline*}
G_{ijk}=\frac{\p}{\p c_i}G_{jk}=\frac{\p}{\p c_i}G_{1,j+k-1}=G_{i,1,j+k-1}=\\=G_{1,i,j+k-1}=\frac{\p}{\p c_1}G_{i,j+k-1}
=\frac{\p}{\p c_1}G_{1,i+j+k-2}=G_{1,1,i+j+k-2}.
\end{multline*}
Similarly we show that all higher derivatives of the form
\[\frac{\p^{r_1}}{\p c_1^{r_1}}\frac{\p^{r_2}}{\p c_2^{r_2}}\dots\frac{\p^{r_n}}{\p c_n^{r_n}}G\]
depend only on $\sum r_i$ and $\sum ir_i$. Since $G$ is a generating function of $g_1,\dots$, the same holds for any individual $g_k$.
The proof of Theorem~\ref{thm:derg} is finished.

\subsection{WDVV-like equation satisfied by the functions $g_i$}

The dependence of $\frac{\partial^3 g_i}{\partial c_k\partial c_l\partial c_m}$ only on the sum $k+l+m$ 
allows us to show that the functions $g_i$, when suitably modified, can satisfy WDVV equation.
Let us pick arbitrary integer $N>2$ and define an $N\times N$ matrix $\eta$ by
\[\eta^{ab}=\begin{cases} 1&\text{ when $a+b=N+1$}\\ 0&\text{ otherwise}\end{cases}\]
\begin{proposition}\label{prop:WDVV}
The function $g_{N+3}$ satisfies the WDVV equation of the form
\[\sum_{\sigma,\tau=1}^N\frac{\partial^3g_{N+3}}{\partial c_\alpha c_\beta c_\sigma}\eta^{\sigma\tau}
\frac{\partial^3g_{N+3}}{\partial c_\mu c_\nu c_\tau}
=
\sum_{\sigma,\tau=1}^N\frac{\partial^3g_{N+3}}{\partial c_\alpha c_\nu c_\sigma}\eta^{\sigma\tau}
\frac{\partial^3g_{N+3}}{\partial c_\mu c_\beta c_\tau}\]
for any $\alpha,\beta,\mu,\nu=1,\dots,N$.
\end{proposition}
\begin{proof}
Let
\[H_s=\frac{\partial^3 g_{N+3}}{\partial c_a\partial c_b\partial c_c},\]
where $a+b+c=s$. We have $H_k=0$ for $k\le 2$ and $k>N+3$.
Let us also define
\[a_1=\alpha+\beta,\,\,a_2=\mu+\nu,\,\,a_3=\alpha+\nu.\]
The statement of Proposition~\ref{prop:WDVV} is trivial when $a_2=a_3$. So let us assume $a_2>a_3$ and $a_1\le a_3$. We need to prove that
\[\sum_{i=1}^{N+1} H_{a_1+i}H_{a_2+N+1-i}=\sum_{i=1}^{N+1}H_{a_3+i}H_{a_1+a_2-a_3+N+1-i}.\]
Substituting $i=j+a_1-a_3$ on the right hand side we get
\[\sum_{i=1}^{N+1}H_{a_1+i}H_{a_2+N+1-i}=\sum_{j=a_3-a_1+1}^{N+a_3-a_1+1}H_{a_1+j}H_{a_2+N+1-j}.\]
Now, for $j>N+1$ we have $a_1+j>N+3$ so $H_{a_1+j}=0$. On the other hand, for $i\le a_3-a_1$ we have $a_2+N+1-i\ge N+1+a_1+a_2-a_3\ge N+a_1+2\ge N+4$,
so $H_{a_2+N-i}=0$. This ends the proof. 
\end{proof}
\begin{remark}
The matrix $\eta^{ab}$ does not come from second derivatives of $g_{N+3}$ as in Section~\ref{sec:reviewWDVV}. 
We could remedy to this problem, by modifying the function
$g_{N+3}$, e.g. adding a new variable. We do not do that, because we do not see any justification for such artificial modifications.
\end{remark}
\begin{remark}
The function $g_{N+3}$ on the space $\C^N$ with basis $c_1,\dots,c_N$ defines a multiplication $\C^N\times\C^N\stackrel{\circ}{\to}\C^N$ 
as explained in Section~\ref{sec:reviewWDVV}. Proposition~\ref{prop:WDVV}
implies then that this multiplication is associative. The homogeneity of $g_{N+3}$ implies that the multiplication is degenerate, which means that there exist 
non--trivial elements $v$
such that the map $\cdot\circ v\colon\C^N\to\C^N$ is not surjective 
(we could take for example $v=c_N$). In the language of WDVV equation and Frobenius manifolds
(compare \cite{Du,He}), such structure is called \emph{massless}.
\end{remark}

\begin{acknowledgement}
The author wishes to express his thanks for D.~Kerner, H.~Żo\-łą\-dek and C.~Hertling for fruitful discussions. 
\end{acknowledgement}

\end{document}